\documentclass[11pt]{article}

\input{amssym.def}
\input{amssym}
\usepackage{eufrak}

\newtheorem{theorem}{Theorem}
\newtheorem{corollary}[theorem]{Corollary}
\newtheorem{definition}[theorem]{Definition}
\newtheorem{lemma}[theorem]{Lemma}
\newtheorem{proposition}[theorem]{Proposition}
\newtheorem{remark}[theorem]{Remark}

\def\M{{\cal{M}}}
\def\N{{\cal{N}}}
\def\P{{\cal{P}}}
\def\tx{\tilde{x}}
\def\u{\mbox{{\sf u}}}
\def\v{\mbox{{\sf v}}}

\def\RRR{\mathbf{R}}
\def\Y{\mathbf{Y}}
\def\YR{\mathbf{Y}(R)}
\def\YYR{\mathbf{Y}_{RRT}(R)}

\def\qq{q^{-1}}

\def\lrh{{\cal L}(R,h)}
\def\lro{{\cal L}(R,1)}
\def\lr{{\cal L}(R)}

\def\De{\Delta}
\def\de{\delta}
\def\la{{\lambda}}

\def\ot{\otimes}
\def\K{{\Bbb K}}
\def\C{{\Bbb C}}
\def\R{{\Bbb R}}

\def\End{{\rm End}}
\def\Ren{R^{\End(V)}}
\def\vv{V^{\otimes 2}}

\def\RR{{\cal R}}

\def\OT{\overline{T}}

\def\be{\begin{equation}}
\def\ee{\end{equation}}

\begin{document}

\makeatletter
\renewcommand{\theequation}{{\thesection}.{\arabic{equation}}}
\@addtoreset{equation}{section} \makeatother

\title{Braided Yangians}

\author{\rule{0pt}{7mm} Dimitri
Gurevich\thanks{gurevich@ihes.fr}\\
{\small\it LAMAV, Universit\'e de Valenciennes,
59313 Valenciennes, France}\\
{\small \it and}\\
{\small \it
Interdisciplinary Scientific Center J.-V.Poncelet}\\
{\small \it Moscow 119002, Russian Federation}\\
\rule{0pt}{7mm} Pavel Saponov\thanks{Pavel.Saponov@ihep.ru}\\
{\small\it
National Research University Higher School of Economics,}\\
{\small\it 20 Myasnitskaya Ulitsa, Moscow 101000, Russian Federation}\\
{\small \it and}\\
{\small \it
Institute for High Energy Physics, NRC "Kurchatov Institute"}\\
{\small \it Protvino 142281, Russian Federation}}

\maketitle

\begin{abstract}
Yangian-like algebras, associated with current $R$-matrices, different from the Yang ones, are introduced. These algebras are of two types. The so-called braided Yangians are close to
the Reflection Equation algebras, arising from involutive or Hecke symmetries. The Yangians of  RTT type  are close to the corresponding RTT algebras. Some properties of these
two classes of the Yangians are studied. Thus, evaluation morphisms
for them are constructed, their bi-algebra structures are described, and quantum analogs of certain symmetric polynomials, in particular, quantum determinants, are introduced.
It is shown that
in any braided Yangian this determinant is always central, whereas in the Yangians of RTT type it is not  in general so.
Analogs of the Cayley-Hamilton-Newton identity in the braided Yangians
are exhibited. A bozonization of the braided Yangians is performed.

\end{abstract}

{\bf AMS Mathematics Subject Classification, 2010:} 81R50, 17B37


\section{Introduction}

The notion of the Yangians was introduced by V.Drinfeld in the framework of his Quantum Group (QG) theory \cite{D}. The Yangian $\Y(gl(N))$ for the general
linear Lie algebra $gl(N)$ is associated with the quantum $R$-matrix, due to Yang
\be
R(u,v)=P-\frac{I}{u-v}.
\label{Ya}
\ee
Hereafter, $I$ is the identity matrix and $P$ stands for the usual flip or its matrix.
The $R$-matrix (\ref{Ya}) meets the quantum Yang-Baxter equation
\be
R_{12}(u,v)R_{23}(u,w)R_{12}(v,w)=R_{23}(v,w)R_{12}(u,w)R_{23}(u,v),
\label{QYBE}
\ee
where $R_{12}=R\ot I$ and $R_{23}=I\ot R$. In general, the lower indices indicate the spaces where a given operator or a matrix is located.

Consider the following system
\be
R(u,v)L_1(u)L_2(v)-L_1(v)L_2(u)R(u,v)=0,
\label{sys}
\ee
where $L(u)=\|l_i^j(u)\|$, $1\leq i, j \leq N$, is a matrix-valued function depending on a formal parameter $u$, $L_1=L\ot I$, $L_2=I\ot L$.

We assume that {\em the generating matrix} $L(u)$ expands in a series
\be
L(u)=\sum_{k=0}^\infty L[k]u^{-k},\quad L[0] = I,
\label{expa}
\ee
where the matrices $L[k]=\|l_i^j[k]\|$ are called the Laurent coefficients of the  matrix $L(u)$.

By one of the equivalent definitions, the Yangian $\Y(gl(N))$ is the unital associative algebra, generated by the elements $l_i^j[k]$, $k\in {\Bbb Z}_+$, subject to the relations
which arise from the system (\ref{sys}). Thus, the current entries $l_i^j(u)$ of the generating matrix $L(u)$ are treated to be elements of $\Y(gl(N))[[u^{-1}]]$. Being presented
in this form, the Yangian is  generated by infinitely many generators $l_i^j[k]$, subject to infinitely many quadratic relations. However, each of these relations contains only a
finite number of generators. In the sequel we deal with other {\em quadratic algebras}, which have the same property.

\begin{remark}\rm
Note that $R$-matrix $\RR(u,v)=PR(u,v)$ is often used instead of $R(u,v)$. Being rewritten in terms of $\RR$, the relation (\ref{sys}) takes the form
\be
\RR(u,v)L_1(u)L_2(v)-L_2(v)L_1(u)\RR(u,v)=0. \label{syss}
\ee
\end{remark}

Now, we list certain basic properties of the Yangians $\Y(gl(N))$.

\begin{enumerate}
\item[{\bf 1.}]
The Yangian has a bi-algebraic structure: the corresponding coproduct is defined on the generating matrix by the rule
$$
\De(L(u))=L(u)\stackrel{.}{\ot} L(u).
$$
Hereafter the notation $A\stackrel{.}{\otimes} B$ stands for the matrix with entries $(A\stackrel{.}{\otimes} B)_i^j=\sum_k\, A_i^k\ot B_k^j $, where
$A$ and $B$ are square matrices of the same size. In terms of the Laurent coefficients the coproduct reads
$$
\Delta(L[m]) = \sum_{k=0}^mL[k]\stackrel{.}{\otimes} L[m-k].
$$
The counit $\varepsilon$ is as follows
$$
\varepsilon(L(u))=I  \quad\Leftrightarrow \quad\varepsilon(l_i^j[m]) = \delta_{m,0}\,\delta_i^j.
$$
\item[{\bf 2.}]
In the Yangian $\Y(gl(N))$ there are well-defined quantum analogs of some symmetric polynomials, namely, the elementary symmetric polynomials and power sums.
The highest quantum elementary symmetric polynomial is called the (quantum) determinant. This determinant (more precisely, its Laurent coefficients)
generate the center of the Yangian $\Y(gl(N))$. The quantum elementary symmetric polynomials generate a commutative subalgebra in the Yangian $\Y(gl(N))$.
The power sums are related with elementary polynomials by a quantum analog of Newton relations and generate the same commutative subalgebra.
\item[{\bf 3.}]
There exists the so-called {\em evaluation} morphism, which sends the Yangian to the enveloping algebra $U(gl(N))$ enabling one to endow any $U(gl(N))$-module with
a $\Y(gl(N))$-module structure.
\end{enumerate}
For these and other properties of the Yangians the reader is referred to \cite{Mo}.

A big interest in the Yangians is motivated in the first turn by their significance in the integrable system theory. In particular, the Yangians play the role of a symmetry group
for the Nonlinear Schr\"{o}dinger hierarchies \cite{MRSZ}. Also, by using the Yangian $\Y(gl(N))$, D.Talalaev \cite{T} succeeded in finding higher Hamiltonians of the
rational Gaudin model.

There are known numerous attempts to generalize the Yangians. Thus, there exist the so-called twisted Yangians which are associated with simple Lie algebras of the series
different from $A_N$ (see \cite{Mo} and the references therein)  and the super-Yangians (see \cite{N,Go}), associated with super-flips $P_{m|n}$.

\begin{remark}\rm
The super-flip  $P_{m|n}$  acts on the tensor square of the super-space $\vv$ by the rule $P_{m|n}(x\ot y)=(-1)^{\overline{x}\,\overline{y}} y\ot  x$. Here $V=V_0\oplus V_1$, $V_0$
(resp., $V_1$) is even (resp., odd) subspace, $\dim V_0=m$, $\dim V_1=n$, elements  $x$ and $y$ are assumed to be homogenous, and $\overline{x}$, $\overline{y}\in\Bbb{Z}_2$
are their parities. We want to emphasize that the super-Yangians from the cited papers are defined via the relation (\ref{syss}) with $\RR=I+\frac{P_{m|n}}{u-v}$. Consequently, they
differ from our super-Yangians of  both types (see below), which are particular cases of our general definitions, provided $R=P_{m|n}$.
\end{remark}

The main purpose  of the current paper is to introduce some new Yangian-like algebras, associated with a large class of $R$-matrices and to study their properties. These algebras are of
two types. Algebras of the first type are similar to the Reflection Equation (RE) algebras\footnote{There are known numerous forms of the RE algebras. In the present paper we use this
term for the algebras associated with constant braidings.} and we call them {\em  braided Yangians}. They are the main objects of the current paper.
Algebras of the second type are similar to  RTT algebras, we  call them {\em Yangians of RTT-type}. These two types of the Yangian-like algebras have different properties. This difference
stems from a similar difference between the RE algebras and RTT ones.

Recall that the RTT and RE algebras are respectively defined by the following systems of relations
\be
R_{12}T_1 T_2-T_1 T_2 R_{12}=0,\qquad T=\|t_i^j\|_{1\leq i,j \leq N}.
\label{RTT}
\ee
\be
R_{12} L_1 R_{12} L_1-L_1 R_{12} L_1 R_{12}=0,\qquad L=\|l_i^j\|_{1\leq i,j \leq N}.
\label{RE}
\ee
Here, $R=R^V:\vv\to \vv,$ is {\it a braiding}. By a braiding we mean a solution of the constant (i.e. without parameters) quantum Yang-Baxter equation
\be
R_{12}R_{23}R_{12} = R_{23}R_{12}R_{23}.
\label{eq:YB}
\ee
and $V$ is a finite-dimensional vector space $V$ over the ground field  $\K$ (we put $\K=\C$ or $\K=\R$),
$\dim_{\K} V=N$,

Mainly, we deal with braidings $R$ obeying  the  condition
\be
(qI-R)(\qq I+R)=0,\qquad q\in \K.
\label{H-cond}
\ee
If $q^2=1$, the corresponding braiding is called  an {\em involutive symmetry}. Unless, it is called {\it a Hecke symmetry}.
Below, if $q\not=1$, $q$ is assumed to be {\em generic}, that is $q^n\not=1$ for any integer $n$.

The best known examples of  Hecke symmetries are standard ones\footnote{Hereafter, we use the term {\em standard} for all objects related to the QG $U_q(sl(N))$.}
However, there exist many other involutive and Hecke symmetries, constructed by one of the authors in the 80's (see \cite{G}).

Also, we assume a basis $\{x_i\}\in V$ to be fixed in the space $V$. Then $R=R^V$ is represented by the corresponding matrix in the basis $\{x_i\ot x_j\}\in \vv$ as follows
$$
R^V(x_i\ot x_j)=R_{ij}^{ks}\,x_k\ot x_s.
$$
The summation over repeated indices is always understood.

The RTT and RE algebras are particular cases of the so-called Quantum Matrix (QM) algebras as defined in \cite{IOP}, but their properties differ drastically. The reader is referred to \cite{GPS}
for a comparative description of these algebras.

The Yangian-like algebras, we are dealing with in the present paper, are current (i.e. depending on a parameter) counterparts of the mentioned QM algebras. In order to define them we have first
to construct a current $R$-matrix $R(u,v)$. By using the so-called {\em Baxterization} procedure, we associate $R(u,v)$ with a given involutive or Hecke symmetry $R$
$$
R(u,v) = R +g(u,v)\,I.
$$
Emphasize that the concrete form of the function $g(u,v)$ depends on $R$ (see Proposition \ref{prop:7}). Now, we are able to give the following definition.

\begin{definition}\rm
Consider a current $N\times N$ matrix $L(u)$ and impose the following  relations
\be
R(u,v)L_1(u) R L_1(v)-L_1(v) R L_1(u) R(u,v)=0.
\label{bY}
\ee
We suppose $L(u)$ to be a series in $u$ analogous to (\ref{expa}). Then, the unital algebra generated by matrix elements $l_i^j[k]$ of the Laurent coefficients $L[k]$, $k\ge 1$,
is called {\it a braided Yangian} and is denoted $\YR$.
\end{definition}

Besides, we introduce an RTT-type Yangian, defined by the system
$$
R(u,v)T_1(u)  T_2(v)-T_1(v)  T_2(u) R(u,v)=0.
$$
As usual, we assume that $T(u)$ expands in a  Laurent series but we do  not impose the condition $T[0]=I$. Such an RTT-type Yangian is denoted $\YYR$. The RTT-type
Yangians corresponding to the standard current $R$-matrices will be called $q$-Yangians.

In this paper we study some basic properties of the braided Yangians and compare them with those of the RTT-type Yangians. First, we define the evaluation morphisms for the
Yangians of both types. It is interesting to note that the concrete form of these morphisms in the braided Yangians depends on the initial symmetry: it differs for the braided Yangians,
arising from involutive symmetries and those arising from Hecke ones. Also, we define the bi-algebra structures in the Yangians of both types. In analogy with the RTT algebra,
the bi-algebra structure of $\YYR$ is usual, while in the braided Yangian $\YR$ it is braided in the sense discussed below.

In all our Yangians we introduce quantum analogs of some symmetric polynomials.
Mainly, we are interested in the elementary symmetric polynomials and, in particular, in the determinant. The quantum determinant, which is the "highest"
elementary symmetric polynomial (it is well-defined if $R$ is even, see the next section), can be defined for a subclass of the Yangians of both types. We show that it is always central in the braided Yangians, whereas in the Yangians of
 RTT type it can be central or not depending on the centrality of the determinant in the corresponding RTT algebra. Also, we exhibit a version of the Cayley-Hamilton-Newton  identity
and the subsequent Cayley-Hamilton and Newton identities for the braided Yangians.

Finally, we consider a bosonization of the braided Yangians. In this connection we discuss  analogs of Zamolodchikov-Faddeev algebras, associated with our current $R$-matrices.

 Note that all properties  of the braided Yangians are more similar to these of the usual Yangians than the corresponding properties of the RTT-type Yangians. Also, note that the
Drinfeld Yangian $\Y(gl(N))$ belongs to the both classes since in this case $R=P$.

The paper is organized as follows. In section 2 we recall some properties of braidings and symmetries, and then we reproduce the main properties of the QM algebras in question (section 3).
In section 4 we exhibit the so-called Baxterization procedure, giving rise to a family of current $R$-matrices. The evaluation morphisms and bi-algebraic structures of all Yangians in
question are discussed in section 5.

In section 6 we define quantum analogs of the matrix powers, skew-matrix powers, some symmetric polynomials,
 and prove the Cayley-Hamilton-Newton identities for the generating matrices of the braided Yangians. Also, as corollaries we get the
 Newton and  Cayley-Hamilton identities.

 In section 7 we perform
a bosonization of the braided Yangians and discuss the Za\-mo\-lod\-chi\-kov-Fad\-deev algebras.

In conclusion of the Introduction we would like to mention some possible applications of the braided Yangians. First, the braided $q$-Yangian gives rise to  interesting
 Poisson structure useful for constructing a new version of the Gaudin model \cite{GS1, GS2}.  Note that the Hamiltonians for this model can be constructed
in the spirit of \cite{T}. Second, by using the Cayley-Hamilton identity for the generating matrices of the braided Yangians it is possible  to perform a quantum analog
of the Drinfeld-Sokolov reduction (see \cite{GST}).

\section{Braidings and symmetries}

In this section we exhibit some properties of symmetries which are used in the sequel. For a more complete list of their properties we refer the reader to the paper \cite{GPS}.
Note that we mainly deal with Hecke symmetries, the corresponding results for involutive ones can be obtained by passing to the limit $q\to 1$.

Below we use a concise notation $R_k:=R_{k\,k+1}$ and the standard notation for $q$-numbers:
$$
k_q=\frac{q^k-q^{-k}}{q-\qq}.
$$
Emphasize, that for generic values of $q$ the $q$-integers are all non-zero: $k_q\not=0$ $\forall\,k\in{\Bbb Z}$.
Now, note  that the following recursive formula
\be
{\cal P}^{(1)} = I,\quad {\cal P}^{(k+1)}_{12\dots k+1} = \frac{k_q}{(k+1)_q}{\cal P}^{(k)}_{12\dots k}\left(\frac{q^k}{k_q} I- R_{k}\right){\cal P}^{(k)}_{12\dots k},\quad k\ge 1
\label{pr-ind}
\ee
defines (for a generic $q$)  some projectors, called {\it skew-symmetrizers}
\be
\P^{(k)}:V^{\ot k}\to  {\bigwedge}_R^{(k)}(V),\quad k\ge 1.
\label{proj}
\ee
Here ${\bigwedge}_R^{(k)}(V)$ is the $k$-th homogenous component of the $R$-skew-symmetric algebra ${\bigwedge}_R(V)$ of the space $V$, where
$$
{\bigwedge}_R(V)=T(V)/\langle {\rm Im}(\qq I+R) \rangle.$$

\begin{definition} \rm We say that a  symmetry $R$ is even if there exists an integer $m\geq 2$ such that
the image of  the skew-symmetrizer $\P^{(m)}$
is non-trivial while the images of the skew-symmetrizers $\P^{(k)}$, $k>m$ are trivial. The couple  $(m|0)$ is called {\it the bi-rank} of this symmetry\footnote{For
a definition of the bi-rank $(m|n)$ the reader is referred to \cite{GPS}.}. \end{definition}

\begin{definition}\rm A braiding $R$ is called {\em skew-invertible},
if there exists an operator $\Psi:\vv\to\vv$ such that
\be
{\rm Tr}_2 R_{12} \Psi_{23}=P_{13}\quad\Leftrightarrow\quad R_{ij}^{kl} \Psi_{lp}^{jq}=\de_i^q \de_p^k.
\label{def:Psi}
\ee
\end{definition}

Given a skew-invertible braiding $R:\vv\to \vv$,  the corresponding   operator $\Psi$ enables us to extend $R$  up to a braiding
\be
R^{V\oplus V^*} :\; (V\oplus V^*)^{\ot 2} \to (V\oplus V^*)^{\ot 2}, \label{exten}
\ee
where $V^*$ is the space dual to $V$, in such a way that the corresponding pairing $V\ot V^*\to \K$ is $R^{V\oplus V^*}$-invariant (see \cite{GPS}).
Emphasize that this extension is unique.  All braidings below are assumed to be skew-invertible.

Also, by means of the  operator $\Psi$ we define two linear operators $B$ and $C$ as follows
\be
B={\rm Tr}_1 \Psi \quad \Leftrightarrow\quad B_i^j=\Psi_{ki}^{kj},\quad {\rm and}\quad  C={\rm Tr}_2 \Psi\quad \Leftrightarrow\quad C_i^j=\Psi_{ik}^{jk}.
\label{BC}
\ee

These operators play an important role in all algebras considered below. In particular, the operator $C$ allows one to introduce the corresponding {\em  $R$-trace} of matrices.
Let $X=\|X_i^j\|$ be an $N\times N$ matrix with entries belonging to any algebra. We define its $ R$-trace as follows
\be
{\rm Tr}_R(X) := {\rm Tr}(X\cdot C) = C_i^j \, X_j^i.
\label{R-tr}
\ee

Let us list some  properties of the operators $B$ and $C$ and of the $R$-trace (see \cite{O, GPS}).  First, for any skew-invertible braiding $R$ one gets
${\rm Tr}_{R(2)} R_{12}=I_1$. Second, for any $N\times N$ matrix $X$ the following holds
 \be
{\rm Tr}_{R(2)}(R^{-1}_{1} X_1 R_{1})= {\rm Tr}_{R(2)}(R_{1}X_1R^{-1}_{1})={\rm Tr}_{R}(X)\, I_1.
\label{basic}
\ee
Third, if a skew-invertible symmetry $R$  is even\footnote{If $R$ is not even,  analogs of formulae (\ref{prR-tr}) and (\ref{prBC}) are presented in \cite{GPS}.}
of bi-rank $(m|0)$, then
\be
{\rm Tr}_R I=q^{-m} m_q.
\label{prR-tr}
\ee
and
\be
B\cdot C=C\cdot B= q^{-2\,m} I.
\label{prBC}
\ee

For a skew-invertible braiding $R$ the corresponding $R$-trace possesses the following cyclic property: for any matrix $X$ of an appropriate size and any polynomial $p(R_1,\dots , R_{k-1})$
we have\footnote{\label{fooo} Observe that formula (\ref{cycl}) is valid even if the matrix $X$ has non-commutative entries. It is so since the matrix $p(R_1,..., R_{k-1})$ is numerical.
Emphasize that the classical formula $Tr\, A\,B=Tr\, B\, A$ is valid if the entries of $A$ commute with these of $B$.}
\be
{\rm Tr}_{R(1...k)} (X\cdot  p(R_1,\dots , R_{k-1}))={\rm Tr}_{R(1...k)}(p(R_1,..., R_{k-1}) \cdot X). \label{cycl} \ee
Hereafter, we use the following  concise notation ${\rm Tr}_{R(1\dots \,k)}: = {\rm Tr}_{R(1)}\dots {\rm Tr}_{R(k)}$.

Indeed, the operators $B$ and $C$ obey the relations (see \cite{O})
\be
R_1B_1B_2 = B_1B_2R_1,\qquad R_1C_1C_2 = C_1C_2R_1.
\label{Ol}
\ee
It is evident that for any positive integer $k$ we have the same type relation
$$
R_kC_kC_{k+1} = C_kC_{k+1}R_k,
$$
and, therefore, the ``string'' $C_1\dots C_n$ commutes with any polynomial in $R_1,\dots R_{n-1}$ and, in particular, with the projectors ${\cal P}^{(k)},\,\, k\leq n-1$.

Now, turn to the role of the operator $B$. Let us endow the space $V^*$ dual to  $V$  with the {\em right dual} basis $\{x^j\}$. Thus, we have a paring
$$
\langle\,,\,\rangle :V\ot V^*\to \K,\quad \langle x_i, x^j\rangle=\de_i^j.
$$
Note that the pairing of these spaces in the reverse order $\langle\,,\rangle: V^*\ot V\to \K$ has to be defined
by $\langle x^j,x_i\rangle=B_i^j$.

We treat elements of the space  $V\ot V^*$ as endomorphisms of the space $V$. To this end we introduce a linear operator by defining the following  action of a basis
element $l_i^j = x_i\otimes x^j$
onto a  vector $x_k\in V$ according to the latter paring:
\be
l_i^j(x_k):= x_i\,\langle x^j, x_k \rangle  =x_i\,  B^j_k.
\label{act}
\ee
The fact that the resulting  map $V\ot V^*\to \End(V)$ is an isomorphism of vector spaces follows from the invertibility of the operator $B$.

Moreover, the usual product (i.e. the composition of endomorphisms)
$$
\circ : \End(V)^{\ot 2}\to \End(V)
$$
is defined in the basis $\{l_i^j\}$ also via the operator $B$:
$$
l_i^j\circ l_k^l=l_i^l \, B^j_k.
$$

If $R$ is a skew-invertible symmetry of bi-rank $(m|0)$, then $\dim \, {\rm Im}\, \P^{(m)}=1$ and consequently
 this skew-symmetrizer can be presented as follows
\be
\P^{(m)}(x_{i_1}\ot\dots \ot x_{i_m}) = \u_{i_1\dots \,i_m} \v,
\label{pr} \ee
where
\be
\v=\v^{j_1\dots j_m}x_{j_1}\ot \dots \ot x_{j_m}\quad {\rm and}\quad \u_{i_1\dots \,i_m}\v^{i_1\dots \,i_m}=1.
\label{vu}
\ee
The tensors $\u$ and $\v$  are defined up to a non-trivial factor: $\u\to \u \la$, $\v\to \v/\la$, $\la\not=0$.
The element
$$
\v=\v^{i_1\dots \,i_m}x_{i_1}\ot\dots \ot x_{i_m}\in {\bigwedge}^{(m)}_R(V)
$$
is a generator of the 1-dimensional  component ${\bigwedge}^{(m)}_R(V)$.

It can be proved that for any $x,y\in V$ there exist $z,u\in V$ such that
$$
R_{m}\dots R_{2}R_{1}(x\ot \v)=\v\ot z,\qquad  R_1R_2\dots R_m(\v\ot y)=u\ot \v.
$$
Therefore, we can define two linear maps $\M,\,\N:V\to V$ by setting $\M(x)=z$, and $\N(y)=u $. As was shown in \cite{G} (in a slightly different normalization) the matrices $\|\M_i^j\|$
and $\|\N_i^j\|$ of these operators are given by the following formulae
\be
\M_i^j=(-1)^{m-1} q\,m_q \,\u_{a_1 a_2\dots a_{m-1} i} \,\v^{j a_1 a_2\dots a_{m-1}},
\label{M}
\ee
\be
\N_i^j=(-1)^{m-1} q\, m_q \,\u_{i a_1 a_2\dots a_{m-1} }\, \v^{a_1 a_2\dots a_{m-1} j}.
\label{N}
\ee

In general, the operators $\M$ and $\N$ are not scalar, but their product is always a scalar operator. Therefore, the  operators $\M$ and $\N$ are scalar simultaneously. The determinant
 in the RTT algebra corresponding to a skew-invertible even symmetry $R$ is central iff these operators are scalar  (see \cite{G} for detail). Below, we prove that it is also true in the corresponding
Yangian of RTT type.

Now, we exhibit  an example of an involutive symmetry $R\in {\rm End}(V^{\otimes 2})$, $\dim V = 2$, for which the operator $\N$ (and consequently $\M$) is not scalar:
\be
R=\left(\begin{array}{cccc}
1&a&-a&ab\\
0&0&1&-b\\
0&1&0&b\\
0&0&0&1
\end{array}
\right), \quad a,b \in \K.
\label{Jord}
\ee

Checking the fact that $R$ is an involutive symmetry is straightforward. Also, it is even and its bi-rank  is $(2|0)$. Moreover, we have
$$
\u_{11}=a,\quad \u_{12}=-1,\quad \u_{21}=1,\quad \u_{22}=0,
\qquad
\v^{11}=0,\quad \v^{12}=-\frac{1}{2},\quad \v^{21}=\frac{1}{2},\quad \v^{22}=-\frac{b}{2}.
$$
The matrix of the operator $\N$ reads
$$
\|\N_i^j\| = \left(\!\!
\begin{array}{cc}
1&a-b\\
0 & 1
\end{array}\!\!\right).
$$
Thus, the  operator $\N$ is scalar iff $a=b$.

This example is very instructive, it shows that the tensor $\v$ does not determine the tensor $\u$ in a unique way (and vice versa). Indeed, if  $b=0$ the symmetric algebra
${\rm Sym}_R(V)=T(V)/\langle \v \rangle$ is classical (it is the polynomial algebra in two commutative generators). However, the tensor $\u$ depends on $a$ and, consequently,
we still have a family of the corresponding involutive symmetries, parameterised by $a$.

Also, note  that the operator $\N$ constructed for a standard Hecke symmetry is scalar. Consequently,  the quantum determinants in the corresponding RTT algebra and
 the $q$-Yangians of RTT type are central.

 The following proposition is important for defining  the determinants in the quantum algebras considered below.

\begin{proposition} \label{profor}
Let $R$ be a skew-invertible even Hecke symmetry of bi-rank $(m|0)$. Then the tensors $\u$ and $\v$
{\rm ({\it see} (\ref{pr}), (\ref{vu}))} satisfy the following relations:
$$
C_{i_1}^{\,j_1}\dots C_{i_m}^{\, j_m}\u_{j_1\dots j_m} = q^{-m^2}\u_{i_1\dots i_m},\qquad C_{i_1}^{\,j_1}\dots C_{i_m}^{\, j_m}\v^{i_1\dots i_m} = q^{-m^2}\v^{j_1\dots j_m},
$$
and consequently
$$
C_1\dots C_m{\cal P}^{(m)} = {\cal P}^{(m)}C_1\dots C_m = q^{-m^2}{\cal P}^{(m)}.
$$
\end{proposition}

\noindent
{\bf Proof.} Since the string $C_1\dots C_m$ commutes  with the projector ${\cal P}^{(m)}$, then using (\ref{pr}) and (\ref{vu}), we get the matrix equality:
$$
C_{i_1}^{\,a_1}\dots C_{i_m}^{\, a_m}\u_{a_1\dots a_m}\v^{j_1\dots j_m} = \u_{i_1\dots i_m}\v^{a_1\dots a_m}C^{\,j_1}_{a_1}\dots C^{\,j_m}_{a_m}.
$$
Contracting the both sides with $\u_{j_1\dots j_m}$ and taking into account (\ref{vu}), we come to
$$
C_{i_1}^{\,a_1}\dots C_{i_m}^{\, a_m}\u_{a_1\dots a_m} = {\rm Tr}_{R(1\dots \,m)}({\cal P}^{(m)})\,\u_{i_1\dots i_m}.
$$
The $R$-trace of ${\cal P}^{(m)}$ can be easily calculated on the base of (\ref{pr-ind}):
$$
{\rm Tr}_{R(m)}({\cal P}^{(m)}_{1\dots m}) = \frac{(m-1)_q}{m_q}\,{\cal P}^{(m-1)}_{1\dots m-1}{\rm Tr}_{R(m)}\Big(\frac{q^{m-1}}{(m-1)_q}I-R_{m-1}\Big)
{\cal P}^{(m-1)}_{1\dots m-1} = \frac{q^{-m}}{m_q}\,{\cal P}^{(m-1)}_{1\dots m-1},
$$
where we used (\ref{prR-tr}) and ${\rm Tr}_{R(m)}R_{m-1} = I_{1\dots m-1}$.

Then, by recurrence we arrive to the following formula
\be
{\rm Tr}_{R(k+1\dots\, m)}({\cal P}^{(m)}_{1\dots m}) = q^{-m(m-k)}\,\frac{k_q!(m-k)_q!}{m_q!}\,{\cal P}^{(k)}_{1\dots k},\qquad 0\le k\le m-1.
\label{part-tr}
\ee
By setting $k=0$ we get the final result
$$
{\rm Tr}_{R(1\dots \,m)}({\cal P}^{(m)}) = q^{-m^2}.
$$
This entails the first claim of the proposition. The second claim can be proved in a similar way.  \hfill\rule{6.5pt}{6.5pt}

\begin{corollary} \label{cor5}
For any matrix $X$ of an appropriate size we have\rm
\be
{\rm Tr}_{R(1\dots \,m)}({\cal P}^{(m)}X) = q^{-m^2}{\rm Tr}_{(1\dots \,m)}({\cal P}^{(m)}X).
\label{R-tr-usual}
\ee
\end{corollary}
Thus, the usual trace and  $R$-trace of the product ${\cal P}^{(m)}X$  differ by a non-trivial factor only.

\section{Quantum matrix algebras related to symmetries}

Before studying the braided Yangians and their RTT counterparts we consider their constant analogs, namely, the RE and RTT algebras. Below, we describe some properties of
these algebras and define quantum analogs of symmetric polynomials.

Let $R$ be any braiding. Together with the RTT and RE algebras consider the unital algebra, generated by elements $l_i^j$, which are subject to the system of relations
\be
R L_1 R L_1-L_1 R L_1 R-h(R L_1-L_1 R)=0, \qquad L=\|l_i^j\|.
\label{mRE}
\ee
This filtered algebra  is called the {\em  modified RE algebra} and is denoted $\lrh$ provided $h\not=0$. If $h=0$, then (\ref{mRE}) turns into defining relations (\ref{RE}), the corresponding
RE algebra will be denoted $\lr$. In what follows, we mainly deal with the algebra $\lro$. The numeric factor $h$ in the right hand side of (\ref{mRE}) enables us to treat  the  algebra $\lro$
as a deformation of $\lr$.

Note that the map
\be
\lr\to \lrh:\qquad L\mapsto L-\frac{h}{q-\qq}I
\label{ma}
\ee
establishes an isomorphism of the algebras $\lr$ and  $\lrh$, provided $q\not=\pm 1$.

Also, note that if $R=R_q$ is a standard Hecke symmetry, then at the limit $q\to 1$ the defining system of the algebra $\lro$ tends to that
$$
P L_1 P L_1-L_1 P L_1 P-(P L_1-L_1 P)=0.
$$
which is nothing but a matrix writing of the relations between the usual generators of the algebra $U(gl(N))$.

Moreover, the algebra $\lrh$ is  covariant with respect to
the adjoint action of the QG $U_q(sl(N))$ (see \cite{GPS}). By contrast, the RTT algebra is not covariant with respect to the adjoint action, it is only covariant with respect to two one-sided actions of
$U_q(sl(N))$.

A very important property of the action  {\rm (\ref{act})} consists in the following: it defines a representation (called {\em covariant}) of the algebra $\lro$ in the space $V$. Whereas the
{\em contravariant} representation of this algebra in the space $V^*$ is defined by the rule
\be
l_i^j(x^l)=-R_{ki}^{lj} \,x^k.
\label{act1}
\ee

Note that the algebra $\lr$ has a braided bi-algebra structure, discovered by Sh.Madjid (see \cite{Maj}). On the generators $l_i^j$ of the algebra $\lr$ the coproduct $\De  : \lr\to \lr\ot \lr$
can be  presented in the following form
\be
\De(L)=L\stackrel{.}{\ot} L\quad \Leftrightarrow \quad \De(l_i^j)=\sum_k l_i^k\ot l_k^j.
\label{cop}
\ee
However, its extension on the higher monomials in generators should be performed via the braiding $\Ren$ (see below (\ref{Ren})).

By using the isomorphism (\ref{ma}) at $h=1$, it is easy to find the coproduct on the generators of the algebra $\lro$:
$$
\De(l_i^j)=l_i^j\ot 1 +1\ot l_i^j-(q-\qq) l_i^k\ot l_k^j.
$$
It also can be extended on the higher monomials in generators by the same braiding $\Ren$.  Consequently, the tensor products $V^{\ot k}\ot {V^*}^{\ot l}$, $\forall\,k,l\in {\Bbb N}$, can be
also endowed with a $\lro$-module structures. According to the representation theory of  the Hecke algebra for a generic $q$ it is possible to decompose the spaces $V^{\ot k}\ot {V^*}^{\ot l}$
into direct sums of some $\lro$-invariant subspaces  in such a way that the final  quasi-tensor rigid  category (called Schur-Weyl one and denoted $SW(V)$) looks like the category of
$U(gl(m|n))$-modules. Besides, with the map (\ref{ma}) we can convert objects  of the category $SW(V)$ into $\lr$-modules. In the sequel the category
$SW(V)$ will also give rise to the corresponding representation categories  of the braided Yangians. Emphasize that  this construction is valid for any skew-invertible Hecke symmetry
$R$ (see \cite{GPS}).

Now, we introduce a useful notation, which allows us to treat  the RE and RTT algebras in a similar manner:
\be
L_{\overline{1}}=L_1,\qquad L_{\overline{k+1}}=R_k L_{\overline{k}}R_k^{-1},\quad k=1,2,...,p-1.
\label{bar-not}
\ee
Note that  $L_{\overline{k}}$ are  $N^p\times N^p$ matrices with entries belonging to  $\lr$ or $\lrh$.

Using this notation we can represent the defining relations of the algebra  $\lrh$  as
\be
R   L_{\overline{1}}L_{\overline{2}} -  L_{\overline{1}}L_{\overline{2}}R-h(L_{\overline{2}}-L_{\overline{1}})=0. \label{mRE1}
\ee
For $h=0$ we get the  relations  in the algebra $\lr$ looking like those in the RTT algebra
$$
R L_{\overline{1}}L_{\overline{2}} -L_{\overline{1}}L_{\overline{2}}R=0.
$$

Let us  exhibit the aforementioned braiding $\Ren$ by using the same notation
\be
\Ren(L_{\overline{1}}\otimes L_{\overline{2}}) = L_{\overline{2}}\otimes L_{\overline{1}}.
\label{Ren}
\ee

As we noticed above, this braiding enters the extension of the coproducts in the algebras $\lr$ and $\lrh$ onto higher monomials.
Thus, for instance, in the algebra $\lr$ we have
\be
\De(L_{\overline{1}}L_{\overline{2}})=(L_{\overline{1}} \ot   L_{\overline{1}}) (L_{\overline{2}} \ot   L_{\overline{2}})=
(L_{\overline{1}}  L_{\overline{2}}) \ot (L_{\overline{1}}  L_{\overline{2}}).
\label{copr}
\ee

Our next aim is to introduce some distinguished elements of the algebra $\lr$, which are quantum analogs of the symmetric polynomials.
\begin{definition}\rm
\label{def:q-elem}
The elements of the algebra $\lr$ defined by
\be
e_k(L) = {\rm Tr}_{R(1\dots\, k)}({\cal P}^{(k)}L_{\overline 1}L_{\overline 2}\dots L_{\overline k})
\label{q-elem}
\ee
are called the {\em (quantum) elementary symmetric polynomials}. If a Hecke symmetry $R$ is even and its  bi-rank  is $(m|0)$, the highest order polynomial $e_m(L)$ is called the {\em quantum determinant} of $L$:
$$
e_m(L) = {\det}_{RE}(L)= {\rm Tr}_{R(1\dots \,m)}(\P^{(m)} L_{\overline{1}}L_{\overline 2}\dots L_{\overline{m}}).
$$
\end{definition}
An important property of $e_k(L)$ is that they belong to the center of the RE algebra $\lr$.

Consider the quantum determinant in more details. In virtue of Corollary \ref{cor5}, with $X=L_{\overline 1}L_{\overline 2}\dots L_{\overline m}$ we can conclude that
$$
{\det}_{RE}(L)=q^{-m^2} {\rm Tr}_{(1\dots \,m)}(\P^{(m)}L_{\overline 1}L_{\overline 2}\dots L_{\overline m}).
$$
Thus, the quantum determinant can be also defined via the usual trace. Moreover, it can be defined without  applying any trace. Namely, we have the following proposition.

\begin{proposition}
\label{prop:77}
The following matrix identity holds true:\rm
\be
{\cal P}^{(m)}L_{\overline 1}L_{\overline 2}\dots L_{\overline m} = L_{\overline 1}L_{\overline 2}\dots L_{\overline m} \,{\cal P}^{(m)} =
q^{m^2}{\det}_{RE}(L)\,{\cal P}^{(m)}.
\label{detL-P}
\ee
\end{proposition}

\noindent
{\bf Proof.} First, we note that the defining relations of the RE algebra $\lr$ can be ``shifted to a higher position":
$$
R_kL_{\overline k}L_{\overline{k+1}} = L_{\overline k}L_{\overline{k+1}}R_k,\quad \forall\,k\ge 1.
$$
Besides, $R_k$ and $L_{\overline m}$ commute with each other if $k\not=m,\, k\not=m+1$.
This entails that  ${\cal P}^{(m)}$ and $L_{\overline 1}L_{\overline 2}\dots L_{\overline m}$ commute:
$$
{\cal P}^{(m)}L_{\overline 1}L_{\overline 2}\dots L_{\overline m} = L_{\overline 1}L_{\overline 2}\dots L_{\overline m}\,{\cal P}^{(m)}.
$$

On taking into account ${\cal P}^{(m)} = {\cal P}^{(m)}{\cal P}^{(m)}$, we get the following
$$
{\cal P}^{(m)}L_{\overline 1}\dots L_{\overline m} = {\cal P}^{(m)} L_{\overline 1}\dots L_{\overline m}{\cal P}^{(m)} =
{\cal P}^{(m)} {\rm Tr}_{(1\dots m)}({\cal P}^{(m)} L_{\overline 1}\dots L_{\overline m}).
$$
It remains to employ the definition of the quantum determinant and formula (\ref{R-tr-usual}).  \hfill\rule{6.5pt}{6.5pt}

\medskip

Denote $H_k(q)$ the Hecke algebra with $k-1$ generators $\sigma_1,\dots ,\sigma_{k-1}$.
This algebra is represented in the spaces $V^{\ot k}$ in a natural way $\rho_R(\sigma_i)=R_i$.

\begin{proposition} {\bf \cite{IP}} Let $R$ be a skew-invertible Hecke symmetry. For an element $z\in H_k(q)$ denote
$$
ch(z) := {\rm Tr}_{R(1\dots \,k)} (Z\,L_{\overline{1}}\dots L_{\overline{k}})={\rm Tr}_{R(1\dots \,k)} (L_{\overline{1}}\dots L_{\overline{k}}\, Z),\quad Z=\rho_R(z).
$$
Consider a linear subspace $Ch_R[L]\subset \lr$ spanned by the unit and the elements $ch(z)$ for all $z\in H_k(q)$, $k\geq 1$.
The space $Ch_R[L]$ is a subalgebra of the center of the algebra $\lr$
\end{proposition}

Observe that the last equality in this formula is valid since the matrix $Z$ is a polynomial in $R_1,\dots, R_{k-1}$. 

The subalgebra  $Ch_R[L]$ is  called  {\em characteristic}.
Note that this claim is also valid for a skew-invertible involutive symmetry. In this case the
role of the algebra $H_k(q)$  is played by  the group algebra of the corresponding symmetric group.

Thus, the elementary symmetric polynomials (\ref{q-elem}) belong to the  characteristic subalgebra of $\lr$. Other distinguished  elements  of  the characteristic subalgebra are
(quantum)  {\em full symmetric polynomials}, and more generally, {\em Schur symmetric polynomials}.

One more family of symmetric polynomials is the so-called (quantum) {\em power sums} defined as follows
 \be
p_k(L)={\rm Tr}_{R(1\dots \,k)}(R_{k-1}\dots R_1\, L_{\overline{1}}\dots L_{\overline{k}})={\rm Tr}_{R(1\dots \,k)}(L_{\overline{1}}\dots L_{\overline{k}}\,R_{k-1}\dots R_1).
\label{pow}
\ee

Emphasize that these power sums can be reduced to the "classical" form $p_k(L)={\rm Tr}_R (L^k)$. Let us demonstrate this property, for instance,  for $k=3$.
Indeed, we have
$$p_3(L)={\rm Tr}_{R(123)}\,L_{\overline{1}}\,R_1 L_{\overline{1}}R^{-1}_1\,R_2  R_1 L_{\overline{1}}R^{-1}_1R^{-1}_2\, R_{2}\, R_1
={\rm Tr}_{R(12)}{\rm Tr}_{R(3)}\,L_{\overline{1}}\,R_1 L_{\overline{1}}R^{-1}_1\,R_2  R_1 L_{\overline{1}}=$$
$${\rm Tr}_{R(12)}\,L_{\overline{1}}\,R_1 L_{\overline{1}}R^{-1}_1\, R_1 L_{\overline{1}}=
{\rm Tr}_{R(12)}\,L_{\overline{1}}\,R_1 L_{\overline{1}} L_{\overline{1}}={\rm Tr}_{R(1)}{\rm Tr}_{R(2)}\, L_{\overline{1}}\,R_1 L_{\overline{1}}^2={\rm Tr}_{R}\,L^3.$$

Now, pass to the    RTT algebras. In these algebras quantum analogs of the above symmetric polynomials can be  defined by similar formulae but
the factors $L_{\overline{k}}$ have to be replaced by $T_k$. Thus, the power sums in the RTT algebras can be defined as follows
\be
p_k(T)={\rm Tr}_{R(1\dots \,k)}(R_{k-1}\dots R_1\, T_{{1}}\dots T_{{k}})={\rm Tr}_{R(1\dots \,k)}(T_{{1}}\dots T_{{k}}\,R_{k-1}\dots R_1).
\label{pow1}
\ee

Also, the corresponding characteristic algebra can be defined in the same way, mutatis mutandis.  The main differences with the RE algebras are the following ones. First,
power sums (\ref{pow1}) cannot be reduced to the form ${\rm Tr}_R\, T^k$ (for $k\geq 2$). Second, the characteristic subalgebra generated by $p_k(T)$ is still commutative
but {\it not} central. It is called {\em a Bethe subalgebra}.

Note that this way of constructing  a Bethe subalgebra in an RTT algebra was suggested by J.-M.Maillet \cite{Ma}. However, in the cited paper  the usual trace was used instead of
its $R$-analog. Note, however, that if in the defining relations (\ref{RTT}) of an RTT algebra one replaces the matrix $T$ by $C\cdot T$, these relations remain valid since the product
$C_1 C_2$ commutes with $R$. So, the map $T\mapsto C\cdot T$ gives rise to an automorphism of the RTT algebra. This entails that though the corresponding Bethe subalgebras differ,
 the proof of their commutativity is the same.

Since the characteristic subalgebra in the RTT algebra is not central, a natural problem arises to describe the center of this algebra. If $R$ is of bi-rank $(m|0)$, the only
non-trivial candidate to the role of a central element is the determinant ${\det}_{RTT}(T)$. However,  as we said above its centrality depends on the
symmetry $R$.

\begin{remark}
\rm
It is possible to define the RTT and RE algebras via the braidings $R$ of the Birman-Wenzl-Murakami type. However, by contrast with the algebras considered above,
if $R$ comes from a QG $U_q(\gg)$ of a series different from $A_N$, the RTT and RE algebras are not deformations of the algebra ${\rm Sym}(\gg)$.
This claim can be proved by checking that the corresponding semiclassical counterparts are not Poisson brackets.
\end{remark}

\section{Baxterization of symmetries}

Now, let us describe the {\em Baxterization} procedure, which associates a current $R$-matrix with a given involutive or Hecke symmetry $R$.
In this section we assume that $\K={\Bbb C}$.

\begin{proposition}\label{prop:7}
Consider the sum\rm
\be
R(u,v)=R+g(u,v)I,
\label{Rmat}
\ee
\it
where $g(u,v)$ is a function of two arguments. Assume that the function $g(u,v)$ depends only on the difference of the arguments: $g(u,v)=f(u-v)$ and that $f(z)$ is a non-constant
meromorphic function. If $R$ is involutive, then $R(u,v)$ is an $R$-matrix iff\rm
\be
g(u,v)=\frac{a}{u-v}.
\label{g}
\ee
\it
If $R$ is a Hecke symmetry, then $R(u,v)$ is an $R$-matrix  iff\rm
\be
 g(u,v)=\frac{q-\qq}{b^{u-v}-1}.
\label{gg}
\ee
\it
Here $a$ and $b\not=1$ are arbitrary nonzero complex numbers.
\end{proposition}

\noindent
{\bf Proof.} Assume $R$ to be a Hecke symmetry. Consider the sum (\ref{Rmat}) and set $g(u, v)=f(u-v)$. Then, imposing the quantum Yang-Baxter equation\footnote{\label{ffo} We use
two notations $R(u,v)$ and $R(x)$ for the $R$-matrices under consideration. Here $x=u-v$ if the function $g(u,v)$ is defined as in (\ref{g}) or (\ref{trig}) and  $x=u/v$  if $g(u,v)$  is of the
form (\ref{rat}). We hope that this notation does not lead to a misunderstanding.}
$$
R_{12}(x)R_{23}(x+y)R_{12}(x)=R_{23}(y)R_{12}(x+y)R_{23}(x),
$$
on  the sum $R(x)=R+f(x)I$ we get the following functional equation on $f(x)$:
$$
(q-q^{-1})f(x+y)+f(x+y)(f(x)+f(y))=f(x)f(y).
$$
Introducing  a new function $h(x)=\frac{q-\qq}{f(x)}$ we get the equation:
\be
h(x)h(y)+h(x)+h(y)=h(x+y).
\label{rel}
\ee

Consequently, the function $\varphi(x)=h(x)+1$ meets the equation $\varphi(x+y)=\varphi(x)\varphi(y)$. This entails $\varphi(x)=b^x$, $b\not=1$.
The case of an involutive symmetry $R$ is left to the reader.
\hfill\rule{6.5pt}{6.5pt}

\medskip

The current $R$-matrix corresponding to (\ref{g}) is called {\em rational}, whereas that corresponding to (\ref{gg}) is called {\em trigonometric}. Note that the Yang $R$-matrices and
their super-analogs are particular cases of  rational braidings. Also, note that the first examples of the Baxterization procedure for general involutive and Hecke symmetries were
presented in \cite{G}.

In general, the parameter $b$ in (\ref{gg}) is independent of $q$, but if we want a trigonometric $R$-matrix to turn into the corresponding rational $R$-matrix as $q\to 1$ (provided a
Hecke symmetry $R_q$ tends to an involutive symmetry when $q\to 1$) we put $b=q^{-2 /a}$. Then we have
\be
R(x)=R_q-\frac{q^{\frac{x}{a}}}{(\frac{x}{a})_q}I,
\label{trig}
\ee
where $x=u-v$. As $q\to 1$ this $R$-matrix tends to
\be
R_1-\frac{a}{x}\, I.
\label{rrat}
\ee
Below, we deal with the current $R$-matrices  (\ref{trig}) and (\ref{rrat}), where for the sake of simplicity we take $a=1$.

Changing  the variables $b^{-u}\to u$, $b^{-v}\to v$ in (\ref{gg}), we get the following form of the trigonometric current $R$-matrix:
\be
R(u,v)=R-\frac{u (q-\qq)}{u-v}I.
\label{rat}
\ee
Evidently, it depends only on the quotient  $x=u/v$.

Consider an example when $R$ is the standard Hecke symmetry with $N=2$. The corresponding current $R$-matrix is (here $x=u/v$):
\be
R(x)=\left(\begin{array}{cccc}
q&0&0&0\\
0&q-\qq&1&0\\
0&1&0&0\\
0&0&0&q
\end{array}\right) -\frac{(q-\qq)x}{x-1}I=
\left(\begin{array}{cccc}
\frac{\qq\, x-q}{x-1}&0&0&0\\
0&\frac{(q-\qq)}{1-x}&1&0\\
0&1&\frac{(q-\qq)x}{1-x}&0\\
0&0&0&\frac{\qq\, x-q}{x-1}
\end{array}\right).
\label{Rm}
\ee

Note that upon multiplying the $R$-matrix $PR(x)$ by an appropriate factor $h(x)$, we can get an $R$-matrix, which unessentially differs from that
from   \cite{FR}.

It is easy to see that the trigonometric current $R$-matrix (\ref{rat}) obeys the following relation
\be
R(u,v)R(v,u)=\Big(1-\frac{uv(q-\qq)^2}{(u-v)^2}\Big) I,
\label{case1}
\ee
whereas the rational current $R$-matrix (\ref{rrat}) with $x=u-v$ is subject to
\be
R(u,v) R(v,u)=\Big(1-\frac{1}{(u-v)^2}\Big) I.
\label{case2}
\ee

Therefore, the following current $R$-matrices
\be
\RRR(u,v)=\frac{R(u,v)}{q-\frac{u(q-\qq)}{(u-v)}}\qquad {\rm and}\qquad \RRR(u,v)=\frac{R(u,v)}{1-\frac{1}{(u-v)}},
\label{demo}
\ee
corresponding to  Hecke symmetries and  involutive ones respectively
turn out to be involutive in the following sense: $\RRR(u,v)\RRR(v,u)=1$.

We use the braidings $\RRR(u,v)$ in the last section in order to define a version of bosonic algebras.

\section{Evaluation morphisms, bi-algebra structures}

Let $R:\vv\to \vv$, $\dim V=N$, be a skew-invertible involutive or Hecke symmetry. Introduce the following notation
$$
[ A, B]_R=R A_{\overline{1}} B_{\overline{2}}-B_{\overline{1}} A_{\overline{2}}R=R A_1R B_1R^{-1}-B_1 R A_1,
$$
where $A$ and $B$ are $N\times N$ matrices. Thus, the expression $[A,B]_R$  is an $N^2\times N^2$ matrix.

For further considerations it is useful to keep in mind the following simple properties valid for an {\it involutive} $R$:
$$
[I,I]_R = [A,I]_R = [I,A]_R = 0,
$$
where $I$ is the unit matrix and $A$ is an arbitrary $N\times N$ matrix. However,  for a {\it Hecke symmetry} the last identity is changed:
$$
[I,I]_R = [A,I]_R = 0,\quad [I,A]_R = (q-q^{-1})(A_{\overline 2} - A_{\overline 1}).
$$

Now, consider the braided Yangian $\YR$, corresponding  to an {\em involutive} symmetry $R$. Taking into account the explicit form of the rational current $R$-matrix $R(u,v)$,
one can write the defining relations of the braided Yangian (\ref{bY}) as follows:
\be
[ L(u), L(v) ]_R =\frac{1}{u-v} \left(L_{\overline{1}}(u)L_{\overline{2}}(v)-L_{\overline{1}}(v)L_{\overline{2}}(u)\right).
\label{new}
\ee

Let us expand the matrix $L(u)$ in a series (\ref{expa}). On substituting this series in (\ref{new}) and equating the coefficients of  $u^{-r} v^{-s}$, we get
$$
\Big[L[r+1],L[s]\Big]_R-\Big[ L[r],L[s+1]\Big]_R =L_{\overline{1}}[r]L_{\overline{2}}[s]-L_{\overline{1}}[s]L_{\overline{2}}[r],\quad
\forall\,r,s\ge 0.
$$
These relations are similar to those from \cite{Mo}. Also, in the same way as in \cite{Mo} we can prove the following claim.

\begin{proposition}  Let $R$ be a skew-invertible involutive symmetry. Then the defining relations of the braided Yangian $\YR$ are equivalent to the system\rm
\be
\Big[ L[r], L[s] \Big]_R =\sum_{a=1}^{{\rm min}(r,\,s)} \left (L_{\overline{1}}[a-1] L_{\overline{2}}[r+s-a]-L_{\overline{1}}[r+s-a] L_{\overline{2}}[a-1]\right),
\quad r,s\ge 1.
\label{def:yang}
\ee

\end{proposition}

The next proposition is important for the representation theory of the braided Yangians.
\begin{proposition}  Let $M$ be the generating matrix of the modified RE algebra  $\lro$ {\rm ({\it see} (\ref{mRE}))} corresponding to a given skew-invertible involutive symmetry $R$.
Then the map\rm
\be
L(u)\mapsto I+\frac{M}{u}
\label{ev}
\ee
\it
defines a surjective morphism $\YR\to \lro$.

Besides, the map $M\mapsto L[1]$ defines an injective morphism $\lro\to \YR$.
\end{proposition}

\noindent
{\bf Proof.}
The defining relations of the algebra $\lro$ are given by (\ref{mRE1}) at $h=1$:
\be
[M,M]_R=M_{\overline{2}}-M_{\overline{1}}.
\label{nnew}
\ee
We have to check that the matrix (\ref{ev}) does satisfy the relation (\ref{new}), provided that (\ref{nnew}) is fulfilled.
This, in turn, is equivalent to checking the relation
$$
(u-v)\left[ \frac{M}{u}, \frac{M}{v} \right]_R=\left(I+\frac{M_{\overline{1}}}{u}\right)\left(I+\frac{M_{\overline{2}}}{v}\right)-
\left(I+\frac{M_{\overline{1}}}{v}\right)\left(I+\frac{M_{\overline{2}}}{u}\right),
$$
which is straightforward.

The second claim of the Proposition is a direct consequence of (\ref{def:yang}) written for $r=s=1$ (recall that $L[0]=I$):
$$
\Big[L[1],L[1]\Big]_R = L_{\overline 2}[1]-L_{\overline 1}[1],
$$
which is equivalent to the defining relations of the algebra $\lro$.
\hfill\rule{6.5pt}{6.5pt}

\medskip

\begin{definition} Similarly to the classical case, the map {\rm (\ref{ev})} is called the evaluation morphism. \end{definition}

Now, let us assume $R$ to be a skew-invertible Hecke symmetry. Taking into account the explicit form (\ref{rat}) of the current $R$-matrix $R(u,v)$,
we have the following defining relations for the braided Yangian $\YR$:
\be
(u-v) [L(u), L(v)]_R=u(q-\qq)\left(L_{\overline{1}}(u)\, L_{\overline{2}}(v)-L_{\overline{1}}(v)\, L_{\overline{2}}(u)\right).
\label{new2}
\ee

By expanding the  matrix $L(u)$ in a series (\ref{expa}), we arrive to the following system
$$
\Big[L[r+1], L[s]\Big]_R-\Big[L[r], L[s+1]\Big]_R= (q-\qq) \left( L_{\overline{1}}[r+1] L_{\overline{2}}[s] -L_{\overline{1}}[s] L_{\overline{2}}[r+1]\right),
\quad r,s\ge 0.
$$

Upon dividing the relation (\ref{new2}) by $u-v$ and using the expansion $\frac{u}{u-v}= \sum_{p=0} \frac{v^p}{u^p}$ we get the following proposition.

\begin{proposition} Let $R$ be a skew-invertible Hecke symmetry. Then the defining relations of the braided Yangian $\YR$ are equivalent to the system
$$
\Big[L[r], L[s]\Big]_R= (q-\qq)  \sum_{a=0}^{{\rm min}(r,\,s-1)} \left( L_{\overline{1}}[a] L_{\overline{2}}[r+s-a] -L_{\overline{1}}[r+s-a] L_{\overline{2}}[a]\right),\quad r, s-1\ge 0.
$$
\end{proposition}

In this case there also exists an evaluation morphism but its surjectivity is not clear, since now the map $M\mapsto L[1]$ does not define a morphism $\lr\to \YR$.

\begin{proposition} Let $M$ be the generating matrix of the  RE algebra  $\lr$ corresponding to a given skew-invertible Hecke  symmetry $R$.
Then the map\rm
\be
 L(u)\mapsto I+\frac{M}{u} \label{ev1}
\ee
\it
defines a  morphism $\YR\to \lr$.
\end{proposition}

\noindent
{\bf Proof.} We have to check that if the matrix $M$ is subject to the relation $[M_{{1}}, M_{{2}}]_R=0$, then
$$
(u-v)\left(R\Big(I+ \frac{M_1}{u}\Big)R\Big(I+ \frac{M_1}{v}\Big)-\Big(I+ \frac{M_1}{v}\Big)R\Big(I+\frac{ M_1}{u}\Big)R\right)=$$
$$(q-\qq)u\, \left( \Big(I+ \frac{M_1}{u}\Big)R\Big(I+ \frac{M_1}{v}\Big)-\Big(I+ \frac{M_1}{v}\Big)R\Big(I+ \frac{M_1}{u}\Big)\right).$$
This equality can be easily verified by a straightforward calculation with the use of the relation $R^2=(q-\qq) R+I$.\hfill\rule{6.5pt}{6.5pt}

\medskip

Let us emphasize the difference between the braided Yangians related to involutive and Hecke symmetries.
In the former case we have the evaluation morphism, which maps the braided Yangian $Y(R)$ into the modified
RE algebra $\lro$, whereas in the latter case it maps the braided Yangian $Y(R)$ into the RE algebra $\lr$.

By using these evaluation morphisms and Schur-Weyl categories of finite dimensional  $\lro$-(or $\lr$-)modules discussed in section 3, we can construct
some representations  of the braided Yangians. Given an object $U$ of this category, we can covert it into a module over the  braided Yangian by putting
\be
l_i^j(u)\triangleright x=x\delta_i^j+\frac{1}{u}\,\rho_U (l_i^j)\triangleright x,\quad x\in U
\label{aact}
\ee
where by $\triangleright$ we denote an action of a linear operator. The symbol $\rho_U$ stands for the corresponding representation of the algebra $\lro$ if $R$ is involutive,
and for that of $\lr$ if $R$ is a Hecke symmetry. Thus, the entries of the matrix $L[1]$ are represented in $U$ according to $\rho_U$, whereas $\rho_U(L[k])=0$ for  $k\geq 2$.

The representation  (\ref{aact})   of the braided Yangian is called the {\em evaluation} representation.

\begin{remark}
\rm If $R$ is involutive, the map $L(u)\to L(u-a),\,\, a\in \K$ is an automorphism of the braided Yangian $\YR$.  Consequently, the  map (\ref{aact}) remains a representation
of the braided Yangian $\YR$, if we replace  the denominator $u$  in  (\ref{aact})  by $u-a$ with  any number $a\in \K$.

If $R$ is a Hecke symmetry and the corresponding braided Yangian $\YR$ is presented in the form (\ref{rat}), then the map $L(u)\to L(au)$, $a\in \K$, $a\not=0$ is an automorphism
of the braided Yangian $\YR$. Consequently, a similar claim is valid for the braided Yangian $\YR$, if we replace  the denominator $u$  in (\ref{aact}) by $au$.

Also, remark that since in general the  braided Yangians do not admit triangular decompositions, analogs of the Drinfeld modules cannot be defined.
\end{remark}

The representations of the braided Yangians described in this section can be multiplied with the use of the following braided bi-algebra structure of $\YR$.
On the generators of $\YR$ the coproduct is defined by:
$$
\De(L(u))=L(u)\stackrel{.}{\ot} L(u).
$$

To extend the above definition on monomials in generators, we again use the operator $\Ren$ discussed in section 3. Thus, similarly to (\ref{copr}) we have
$$
\De(L_{\overline{1}}(u) L_{\overline{2}}(v))=(L_{\overline{1}}(u) \ot L_{\overline{1}}(u)) (L_{\overline{2}}(v) \ot L_{\overline{2}}(v))=
(L_{\overline{1}}(u)  L_{\overline{2}}(v))\ot (L_{\overline{1}}(u) L_{\overline{2}}(v))
$$
and so on. Thus, the coproduct is similar to that in the RE algebra $\lr$.

As for the RTT-type Yangians, the corresponding coproduct is similar to  that in the RTT algebras and the bi-algebra structure is usual.

However, in order to define  the evaluation morphisms for these Yangians we  do not impose the condition $T[0]=I$ in the expansion of the
function $T(u)$. Then we introduce the evaluation map by
$$
T(u)\to T+\frac{\OT}{u}.
$$
It is a morphism of algebras if the target algebra is defined by the relations
$$
R\,T_1\, T_2=T_1\, T_2\, R,\qquad R\,\OT_1\, \OT_2=\OT_1\, \OT_2\, R,\qquad R\,\OT_1\, T_2=T_1\, \OT_2\, R,
$$
provided $R$ is a Hecke symmetry, and
$$
R\,T_1\, T_2=T_1\, T_2\, R,\qquad R\,\OT_1\, \OT_2-\OT_1\, \OT_2\, R=T_1\, \OT_2-\OT_1\, T_2,\qquad R\,\OT_1\, T_2=T_1\, \OT_2\, R,
$$
provided $R$ is involutive. In the standard case the former algebra is similar to the QG $U_q(gl(m))$, as presented in \cite{Mo}. We do not know any treatment of the latter algebra.

\section{Determinant, symmetric polynomials, CHN identities}

In this section we define quantum analogs of some symmetric polynomials in the braided Yangians $\YR$ and show that the quantum determinant of the generating matrix $L(u)$
is always central in $\YR$. Recall that the quantum determinant is the highest order elementary symmetric polynomial, which is well-defined provided  the initial  Hecke
symmetry $R$ is  even. Also, we establish the so-called Cayley-Hamilton-Newton identities for the matrix $L(u)$ and their corollaries --- the Newton relations and the Cayley-Hamilton identity.

Assume $R$ to be current $R$-matrix associated with a Hecke symmetry. Then according to footnote \ref{ffo} we present the defining relations of the corresponding Yangian $\Y(R)$ as follows
\be
R(u/v)L_{\overline 1}(u)L_{\overline 2}(v) = L_{\overline 1}(v)L_{\overline 2}(u)R(u/v),\qquad R(x)=R-\frac{(q-\qq)x}{x-1}.
\label{Br-Y-ndef}
\ee
Note, that setting $x=q^{2k}$, we get $R(q^{2k}) = R-\frac{q^k}{k_q}\,I$.

Now, we exhibit some technical formulae which will be used below. They  can be easily proved by induction.
\begin{eqnarray}
{\cal P}^{(k+1)}_{1\dots k+1}&=&\frac{(-1)^k}{(k+1)_q}\,R_1(q^2)R_2(q^4)\dots R_k(q^{2k})\,{\cal P}^{(k)}_{12\dots k}\nonumber\\
&=& \frac{(-1)^k}{(k+1)_q}\,R_k(q^2)R_{k-1}(q^4)\dots R_1(q^{2k})\,{\cal P}^{(k)}_{2\dots k+1}\label{q-anti}\\
&=& \frac{(-1)^k}{(k+1)_q}\,{\cal P}^{(k)}_{12\dots k}\,R_k(q^{2k})R_{k-1}(q^{2(k-1)})\dots R_1(q^2)\nonumber\\
&=& \frac{(-1)^k}{(k+1)_q}\,{\cal P}^{(k)}_{2\dots k+1}\,R_1(q^{2k})R_2(q^{2(k-1)})\dots R_k(q^2).\nonumber
\end{eqnarray}

To simplify the writing we introduce the following compact notation
\be
R^{(\pm)}_{i\rightarrow j}(u) =
\left\{
\begin{array}{l}
R_i(u)R_{i+1}(q^{\pm 2}u)\dots R_j(q^{\pm 2(j-i)}u),\quad {\rm if}\quad j\ge i\\
\rule{0pt}{7mm}
R_i(u)R_{i-1}(q^{\pm 2}u)\dots R_j(q^{\pm 2(i-j)}u),\quad {\rm if}\quad i>j.\\
\end{array}\right.
\label{R-string}
\ee
For example, in this notation the first line in (\ref{q-anti}) can be presented as follows
$$
{\cal P}^{(k+1)}_{1\dots k+1} = \frac{(-1)^k}{(k+1)_q}\,R^{(+)}_{1\rightarrow k}(q^2)\,{\cal P}^{(k)}_{12\dots k}.
$$

With the use of (\ref{q-anti}) we can prove the following lemma.

\begin{lemma}\label{lem:17}
For any value of the parameter $u$ and for any integer $k\ge 1$ the following holds
$$
R_{1\rightarrow k}^{(+)}(q^{-2(k-1)}u)\,{\cal P}^{(k)}_{1\dots k } =
{\cal P}^{(k)}_{2\dots k+1 }\,R_{1\rightarrow k}^{(-)}(u).
$$
If the rank of a Hecke symmetry is $(m|0)$, then for the highest order skew-symmetrizer ${\cal P}^{(m)}$ the following identities hold true:\rm
\be
\begin{array}{rcr}
R_{1\rightarrow m}^{(+)}(q^{-2(m-1)}u)\,{\cal P}^{(m)}_{1\dots m }\!\!\!
&=&\!\!\!\displaystyle
(-1)^{m-1}q\,m_q\, \frac{(u-q^{2m})}{(q^2 u-q^{2m})}\,{\cal P}^{(m)}_{2\dots m+1}{\cal P}^{(m)}_{1\dots m},\\
\rule{0pt}{8mm}
{\cal P}^{(m)}_{2\dots m+1 }R_{1\rightarrow m}^{(-)}(u)\!\!\! &=&\!\!\! \displaystyle
(-1)^{m-1}q\,m_q\, \frac{(u-q^{2m})}{(q^2 u-q^{2m})}\,{\cal P}^{(m)}_{2\dots m+1}{\cal P}^{(m)}_{1\dots m}.
\end{array}
\label{q-anti-m}
\ee
\end{lemma}

Now, we define analogs of the elementary symmetric polynomials for the braided Yangian $\YR$ associated with a skew-invertible Hecke symmetry of bi-rank $(m|n)$.
\begin{definition}\rm
\label{def:Y-elem}
Let $L(u)$ be the generating matrix of a braided Yangian, associated with a Hecke symmetry $R$. The elements
\be
e_0(u) = 1,\quad e_k(u) = {\rm Tr}_{R(1\dots k)}\Big({\cal P}^{(k)}_{1\dots k}\,L_{\overline 1}(u)L_{\overline 2}(q^{-2}u)\dots L_{\overline k}(q^{-2(k-1)}u)\Big),\quad k\geq 1
\label{elem-sym}
\ee
will be called the {\em elementary symmetric polynomials}.

If the rank of a Hecke symmetry is $(m|0)$, then the highest order nonzero elementary symmetric polynomial $e_m(u)$ is called the {\em determinant} of $L(u)$ and
denoted $\det_{\YR}(u)$.
\end{definition}
An important property of the elements $e_k(u)$ is that they commute with each other and consequently their Laurent coefficients in the corresponding series in $u^{-1}$
generate a commutative subalgebra in the braided Yangian $\YR$.
The proof of this fact, as well as the analysis of the structure of the commutative subalgebra, will be given in our subsequent publication.

Now, we consider the highest order elementary symmetric polynomial $e_m(u) =\det_{\YR}(u)$ and prove its main property.
\begin{proposition}
If $R$ is a skew-invertible even Hecke symmetry of bi-rank $(m|0)$ then in the corresponding braided Yangian $\Y(R)$ the determinant $e_m(u) = \det_{\YR}(u)$
defined by (\ref{elem-sym}) with $k=m$ is central, i.e.
\rm
\be
{\det}_{\YR}(u)L(v) = L(v){\det}_{\YR}(u), \quad \forall\, u,v.
\label{det-cent}
\ee
\end{proposition}

\noindent
{\bf Proof.}
The braided Yangian defining relations (\ref{Br-Y-ndef}) allow one to write the following identity which is valid for arbitrary values of the parameters $u$ and $v$:
\begin{eqnarray*}
R_{1\rightarrow m}^{(-)}(u/v)L_{\overline 1}(u)L_{\overline 2}(q^{-2}u)&\!\!\!\!\dots\!\!\!\!& L_{\overline m}(q^{-2(m-1)}u)L_{\overline{m+1}}(v)\\
&\!\!\!\!=\!\!\!\!&L_{\overline 1}(v)L_{\overline 2}(u)L_{\overline 3}(q^{-2}u)\dots L_{\overline{ m+1}}(q^{-2(m-1)}u)R_{1\rightarrow m}^{(-)}(u/v).
\end{eqnarray*}

We multiply the both sides of this identity by ${\cal P}^{(m)}_{2\dots m+1}$ on the left hand side and then use the second relation of (\ref{q-anti-m}):
\begin{eqnarray}
\phi(u/v){\cal P}^{(m)}_{2\dots m+1}{\cal P}^{(m)}_{1\dots m}L_{\overline 1}(u)\!\!\!\!&\dots& \!\!\!\!L_{\overline m}(q^{-2(m-1)}u)L_{\overline {m+1}}(v)\label{prom}\\
&=&\!\!\!\!L_{\overline 1}(v){\cal P}^{(m)}_{2\dots m+1}L_{\overline 2}(u)\dots L_{\overline {m+1}}(q^{-2(m-1)}u)R_{1\rightarrow m}^{(-)}(u/v),\nonumber
\end{eqnarray}
where
$$
\phi(u) = (-1)^{m-1}q\,m_q\, \frac{(u -q^{2m})}{(q^2 u-q^{2m})}
$$
stands for the scalar factor in the right hand side of (\ref{q-anti-m}).

Taking into account the defining relations (\ref{Br-Y-ndef}) and the recurrent formulae (\ref{q-anti}) we get
\be
{\cal P}^{(m)}_{1\dots \,m}L_{\overline 1}(u)L_{\overline 2}(q^{-2}u)\dots L_{\overline m}(q^{-2(m-1)}u) = L_{\overline 1}(q^{-2(m-1)}u)\dots
L_{\overline{m-1}}(q^{-2}u)L_{\overline m}(u)\,{\cal P}^{(m)}_{1\dots \,m}.
\label{PL-str}
\ee
Then in full analogy with Proposition \ref{prop:77} we conclude that
\be
{\cal P}^{(m)}_{1\dots m}L_{\overline 1}(u)\dots L_{\overline m}(q^{-2(m-1)}u) = q^{m^2}e_m(u){\cal P}^{(m)}_{1\dots m}.
\label{Pm-det}
\ee
So, the left hand side of (\ref{prom}) takes the form:
$$
\phi(u/v)q^{m^2}{\cal P}^{(m)}_{2\dots m+1}{\cal P}^{(m)}_{1\dots m}e_m(u)L_{\overline {m+1}}(v).
$$

Transforming the right hand side of (\ref{prom}) in the same way we obtain:
$$
L_{\overline 1}(v){\cal P}^{(m)}_{2\dots m+1}L_{\overline 2}(u)\dots L_{\overline {m+1}}(q^{-2(m-1)}u)R_{1\rightarrow m}^{(-)}(u/v) =
\phi(u/v)q^{m^2}L_{\overline 1}(v)e_m(u){\cal P}^{(m)}_{2\dots m+1}{\cal P}^{(m)}_{1\dots m}.
$$
Therefore, we come to the following intermediate result:
\be
{\cal P}^{(m)}_{2\dots m+1}{\cal P}^{(m)}_{1\dots m}e_m(u)L_{\overline {m+1}}(v) =
L_{\overline 1}(v)e_m(u){\cal P}^{(m)}_{2\dots m+1}{\cal P}^{(m)}_{1\dots m}.
\label{Le}
\ee

Now, on taking into account that $L_{\overline{m+1}}(v) = R_{m}\dots R_1L_1(v)R^{-1}_1\dots R^{-1}_m$, we multiply the identity (\ref{Le})
by $R_m\dots R_1$ from the right and use the following properties of the highest skew-symmetrizer (see \cite{H-T})
$$
{\cal P}^{(m)}_{1\dots m}R_m\dots R_1 = (-1)^{m-1}q\,m_q{\cal P}^{(m)}_{1\dots m}{\cal P}^{(m)}_{2\dots m+1},
$$
$$
{\cal P}^{(m)}_{2\dots m+1}{\cal P}^{(m)}_{1\dots m}{\cal P}^{(m)}_{2\dots m+1} = m_q^2\,{\cal P}^{(m)}_{2\dots m+1}.
$$
Finally, we get
$$
{\cal P}^{(m)}_{2\dots m+1}e_m(u)L_1(v) = L_1(v)e_m(u){\cal P}^{(m)}_{2\dots m+1}.
$$

Since ${\rm Tr}_{R(2\dots m+1)}{\cal P}^{(m)}_{2\dots m+1} = q^{-m^2}$ is nonzero, then on applying  the corresponding $R$-trace to the both sides of the above
equality we come to the matrix identity
$$
e_m(u)L_1(v) = L_1(v)e_m(u).
$$
This completes the proof.\hfill \rule{6.5pt}{6.5pt}

\medskip

\begin{remark}\rm
Note, that the same considerations can be carried out for any RTT-type Yangian with the defining relations
$$
R(u/v)T_1(u)T_2(v) = T_1(v)T_2(u)R(u/v).
$$
In the same way we can come to the intermediate formula analogous to (\ref{Le}):
$$
{\cal P}^{(m)}_{2\dots m+1}{\cal P}^{(m)}_{1\dots m}e_m(u)T_{m+1}(v) =
T_1(v)e_m(u){\cal P}^{(m)}_{2\dots m+1}{\cal P}^{(m)}_{1\dots m},
$$
where the elementary symmetric polynomials $e_k(u)$ in this case are defined by the relation
\be
e_k(u) = {\rm Tr}_{R(1\dots k)}\Big({\cal P}^{(k)}_{1\dots k}\,T_1(u)T_2(q^{-2}u)\dots T_k(q^{-2(k-1)}u)\Big).
\label{RTT-elem}
\ee

Writing the skew-symmetrizer ${\cal P}^{(m)}$ in terms of structural tensors $\u$ and $\v$ (\ref{pr}) and using the definition (\ref{N}) of the matrix $\cal N$
we can transform the above identity to the following form:
$$
{\cal N}\cdot T(v)e_m(u) = e_m(u)T(v)\cdot {\cal N}.
$$
So, in the RTT-type Yangian the determinant $e_m(u)$ is central iff the matrix $\cal N$ (and, therefore, $\cal M$) is scalar.
\end{remark}

Now, we go back to the Yangians $\Y(R)$. Let us introduce the skew-symmetric powers of the generating matrix $L(u)$ by the rule:
\be
L^{\wedge 1}(u) = L(u),\quad
L^{\wedge k}(u) = {\rm Tr}_{R(2\dots k)}\Big({\cal P}^{(k)}_{12\dots k}L_{\overline 1}(u)L_{\overline 2}(q^{-2}u)\dots L_{\overline k}(q^{-2(k-1)}u)\Big)
\quad k\ge 2.
\label{wedge-p}
\ee
If the bi-rank of the Hecke symmetry $R$ is $(m|0)$, then $L^{\wedge k}(u)$ is equal to zero for $k>m$. Also, it is evident
that $e_k(u) = {\rm Tr}_R(L^{\wedge k}(u))$  for any $k\ge 1$.

Besides, we define quantum analogs of the matrix powers of the generating matrix:
\be
L^{[k]}(u) = {\rm Tr}_{R(2\dots k)}\Big(L_{\overline 1}(q^{-2(k-1)}u)L_{\overline 2}(q^{-2(k-2)}u)\dots L_{\overline k}(u)R_{k-1}\dots R_2R_1\Big),
\quad k\ge 1.
\label{power-p}
\ee
Taking into account the definition of $L_{\overline s}$ and the property ${\rm Tr}_{R(s+1)}R_s = I$, the above expression for the matrix power
 can be reduced to the product
\be
L^{[k]}(u) = L(q^{-2(k-1)}u)L(q^{-2(k-2)}u)\dots L(u).
\label{mult-pow}
\ee
This reducing is realised in the same way as it was done in Section 3 while computing $p_3(L)$.
 However, in the braided Yangians the  {\em quantum matrix powers} $L^{[k]}(u)$ differ from the classical ones
$L^{k}(u)$ by  multiplicative shifts of the parameters in the matrices $L(u)$.

By analogy with the elementary symmetric polynomials we define {\it the power sums}:
\be
p_k(u) = {\rm Tr }_R(L^{[k]}(u)).
\label{pow-sum}
\ee
\begin{proposition}\label{prop:22}
The skew-symmetric and matrix powers of $L(u)$ satisfies the series of Cayley-Hamilton-Newton identities \rm
\be
(-1)^{k+1}k_qL^{\wedge k}(u) = \sum_{p=0}^{k-1}(-q)^pL^{[k-p]}(q^{-2p}u)e_p(u),\quad k\ge 1.
\label{chn-id}
\ee
\end{proposition}

\noindent{\bf Proof.}
First of all, taking into account (\ref{PL-str}), we can identically rewrite the definition (\ref{wedge-p}) in the form:
\be
L^{\wedge k}(u) = {\rm Tr}_{R(2\dots k)}\Big(L_{\overline 1}(q^{-2(k-1)}u)L_{\overline 2}(q^{-2(k-2)}u)\dots L_{\overline k}(u){\cal P}^{(k)}_{12\dots k}\Big).
\label{new-w-p}
\ee
Then we use the recurrent formula for the skew-symmetrizers ${\cal P}^{(k)}$ equivalent to (\ref{pr-ind}):
$$
k_q{\cal P}^{(k)}_{1\dots k} = {\cal P}^{(k-1)}_{2\dots k}\Big(q^{k-1} - (k-1)_qR_1\Big){\cal P}^{(k-1)}_{2\dots k}.
$$
So, on aplying the above relation to (\ref{new-w-p}) we obtain:
\begin{eqnarray}
k_qL^{\wedge k}(u)\!\!&=&\!\!\! q^{k-1}L_{\overline 1}(q^{-2(k-1)}u) e_{k-1}(u)\label{L-w}\\
&-&\!\!\!(k-1)_q{\rm Tr}_{R(2\,\dots\, k)}\Big(L_{\overline 1}(q^{-2(k-1)}u) L_{\overline 2}(q^{-2(k-2)}u)\dots L_{\overline k}(u)
{\cal P}^{(k-1)}_{2\dots k}R_1{\cal P}^{(k-1)}_{2\dots k}\Big).\nonumber
\end{eqnarray}
Here, to get the first term in the right hand side we applied the identity
$$
{\rm Tr}_{R(r\dots\, r+k)}\Big({\cal P}^{(k+1)}_{r\dots\, r+k}L_{\overline r}(u)\dots L_{\overline {r+k}}(x+k)\Big) = e_{k+1}(u)I_{1\dots \,r-1}, \quad \forall \, r\ge 1.
$$

Let us simplify the second term in the right hand side of (\ref{L-w}). For this purpose we use the cyclic property of the $R$-trace:
$$
{\rm Tr}_{R(r\dots\, r+k)}\Big({\cal P}^{(k+1)}_{r\dots \,r+k}\,X\Big) = {\rm Tr}_{R(r\dots \,r+k)}\Big(X{\cal P}^{(k+1)}_{r\dots \,r+k}\Big),
$$
where $X$ stands for any matrix of an appropriate size. Consequently, we have
\begin{eqnarray*}
&{\rm Tr}_{R(2\dots k)}&\!\!\!\!\!\!\Bigr(
{\cal P}^{(k-1)}_{2\dots k}L_{\overline 1}(q^{-2(k-1)}u)L_{\overline 2}(q^{-2(k-2)}u)\dots L_{\overline k}(u){\cal P}^{(k-1)}_{2\dots k}R_1\Big) =(\mbox{\rm applying (\ref{PL-str})})\\
&&={\rm Tr}_{R(2\dots k)}\Big(
{\cal P}^{(k-1)}_{2\dots k}L_{\overline 1}(q^{-2(k-1)}u)L_{\overline 2}(u)\dots L_{\overline k}(q^{-2(k-2)}u)R_1\Big) = (\mbox{\rm applying (\ref{PL-str}) again})\\
&&={\rm Tr}_{R(2\dots k)}\Big(
L_{\overline 1}(q^{-2(k-1)}u)L_{\overline 2}(q^{-2(k-2)}u)\dots L_{\overline k}(u){\cal P}^{(k-1)}_{2\dots k}R_1\Big).
\end{eqnarray*}
Besides, in the second line we have taken into account that ${\cal P}^{(s)}{\cal P}^{(s)} = {\cal P}^{(s)}$ for any $s$.

So, we come to the following intermediate result:
\begin{eqnarray*}
k_qL^{\wedge k}(u)\!\!&=&\!\! q^{k-1}L_{\overline 1}(q^{-2(k-1)}u) e_{k-1}(u)\\
&-&\!\!(k-1)_q{\rm Tr}_{R(2\dots k)}\Big(
L_{\overline 1}(q^{-2(k-1)}u)L_{\overline 2}(q^{-2(k-2)}u)\dots L_{\overline k}(u){\cal P}^{(k-1)}_{2\dots k}R_1\Big).
\end{eqnarray*}

Now, we apply the same transformations to $(k-1)_q{\cal P}^{(k-1)}_{2\dots k}$ and so on. The required result (\ref{chn-id}) is proved by induction
with the use of definitions of $e_k(u)$ and $L^{[k]}(u)$.\hfill \rule{6.5pt}{6.5pt}

\medskip

\begin{remark}
\rm
Note, that formula (\ref{chn-id}) is valid for {\it any} positive integer $k$ and for any skew-invertible Hecke symmetry $R$. However, if the bi-rank of a Hecke
symmetry is $(m|0)$, then we have $L^{\wedge k}\equiv 0$ for $k>m$, and all matrix identities given by (\ref{chn-id}) at $k>m$ are consequences of the
Cayley-Hamilton identity corresponding to $k=m$ (see Proposition \ref{prop:CH-id} below).
\end{remark}
As a simple corollary of the Proposition \ref{prop:22} we get analogs of the Newton identities, connecting the elementary symmetric polynomials and
power sums.
\begin{corollary}\rm
The elementary symmetric polynomials (\ref{elem-sym}) and power sums (\ref{pow-sum}) are connected by the set of identities ($k\ge 1$):
\be
k_qe_k(u) - q^{k-1}p_1(q^{-2(k-1)}u)e_{k-1}(u)+q^{k-2}p_2(q^{-2(k-2)}u)e_{k-2}(u)+\dots +(-1)^kp_k(u)= 0.
\label{q-Newt}
\ee
\end{corollary}
These relations can be easily proved by taking the $R$-trace of the identities
(\ref{chn-id}).


Now, consider the  Cayley-Hamilton-Newton identity (\ref{chn-id}) at $k=m$. Actually, it only  contains the matrix powers
$L^{[p]}$ and is called the Caley-Hamilton identity for the matrix $L(u)$.
\begin{proposition}\label{prop:CH-id}
The generating matrix of the braided Yangian $\YR$ associated with the Hecke symmetry $R$ of bi-rank $(m|0)$ meets the matrix Cayley-Hamilton
identity:\rm
\be
\sum_{p=0}^m(-q)^{p}L^{[m-p]}(q^{-2p}u)e_p(u) = 0,
\label{CayHam}
\ee\it
where by definition we set $L^{[0]}(u) := I$.
\end{proposition}

\noindent
{\bf Proof.} Let us set $k=m$ in (\ref{chn-id}):
\be
(-1)^m m_qL^{\wedge m}(u)+\sum_{p=0}^{m-1}(-q)^{p}L^{[m-p]}(q^{-2p}u)e_p(u) = 0.
\label{chn-m}
\ee
Transform now the first term of this identity, using the relation (\ref{Pm-det}) valid for the highest order skew-symmetrizer:
\begin{eqnarray*}
m_q L^{\wedge m}(u)\!\!\!\! &=& \!\!\!\! m_q{\rm Tr}_{R(2\,\dots \,m)}\Big({\cal P}^{(m)}_{12\,\dots\, m}L_{\overline 1}(u)L_{\overline 2}(q^{-2}u)\dots L_{\overline m}(q^{-2(m-1)}u)\Big)\\
\!\!\!\! &\stackrel{(\ref{Pm-det})}{=}& \!\!\!\! q^{m^2}m_qe_m(u){\rm Tr}_{R(2\,\dots \,m)}\Big({\cal P}^{(m)}_{12\,\dots\, m}\Big) = q^me_m(u)I_1.
\end{eqnarray*}

Here, at the last step we use the value of a multiple trace of the skew-symmetrizer
$$
m_q{\rm Tr}_{R(2\,\dots \,m)}\Big({\cal P}^{(m)}_{12\,\dots\, m}\Big) = q^{-m(m-1)}I_1,
$$
which is a particular case of (\ref{part-tr}).

Substituting the relation $m_q L^{\wedge m}(u) = q^me_m(u)I$ into (\ref{chn-m}) we get the identity (\ref{CayHam}).\hfill\rule{6.5pt}{6.5pt}

\medskip

\begin{remark}\rm
For the RTT-type Yangian the Cayley-Hamiltion identity slightly differs in the last term, corresponding to the zero power of $T(u)$. Indeed, for this type Yangian
we have
$$
m_qT^{\wedge m}(u) = q^m e_m(u){\cal N},
$$
and the CH identity in this case reads
$$
q^m e_m(u){\cal N}+\sum_{p=0}^{m-1}(-q)^{p}T^{[m-p]}(q^{-2p}u)e_p(u) = 0,
$$
where the elements $e_k(u)$ are defined by (\ref{RTT-elem}) and $T^{[k]}$ by
$$
T^{[k]}(u) = {\rm Tr}_{R(2\dots k)}\Big(T_1(q^{-2(k-1)}u)T_2(q^{-2(k-2)}u)\dots T_k(u)R_{k-1}R_{k-2}\dots R_1\Big).
$$
\end{remark}

The case when $R$ is involutive and consequently the corresponding current $R$-matrix is rational can be treated in a similar manner.
Let us define the matrix skew-powers $L^{\wedge k}(u)$ and matrix powers $L^{[k]}(u)$  by the formulae
analogous to (\ref{wedge-p}) and (\ref{power-p}):
\be
L^{\wedge 1}(u) = L(u),\quad
L^{\wedge k}(u) = {\rm Tr}_{R(2\dots k)}\Big({\cal P}^{(k)}_{12\dots k}L_{\overline 1}(u)L_{\overline 2}(u-1)\dots L_{\overline k}(u-k+1)\Big)
\quad k\ge 2,
\label{inv-wedge-p}
\ee
\be
L^{[k]}(u) = {\rm Tr}_{R(2\dots k)}\Big(L_{\overline 1}(u-k+1)L_{\overline 2}(u-k+2)\dots L_{\overline k}(u)R_{k-1}R_{k-2}\dots R_1 \Big).
\quad k\ge 1.
\label{inv-power-p}
\ee
Here the projectors ${\cal P}^{(k)}$ are obtained as the limit $q\rightarrow 1$ of (\ref{pr-ind}). The corresponding elementary symmetric
polynomials and power sums are respectively defined by  $e_k(u) = {\rm Tr}_RL^{\wedge k}(u)$ and $p_k(u) = {\rm Tr}_RL^{[k]}(u)$.

Note that in this case the {\em quantum matrix powers} $L^{[k]}(u)$ can be reduced to a form similar to (\ref{mult-pow}) but with additive shifts of the parameter.

\begin{proposition}
\label{prop:inv-ident}
In the braided Yangian $\Y(R)$ associated with an involutive skew-invertible symmetry $R$ of bi-rank $(m|0)$ the following
set of Cayley-Hamilton-Newton relations exists:\rm
\be
(-1)^{k+1}kL^{\wedge k}(u) = \sum_{p=0}^{k-1}(-1)^{p}L^{[k-p]}(u-p)e_p(u),\quad k\ge 1.
\label{inv-chn-id}
\ee
\it At $k=m$ we get the Cayley-Hamilton identity
\be
\sum_{p=0}^m(-1)^{p}L^{[m-p]}(u-p)e_p(u) = 0.
\label{CayHam1}
\ee
The elementary symmetric polynomials and power sums are connected by the Newton relations
\be
ke_k(u) - p_1(u-k+1)e_{k-1}(u)+p_2(u-k+2)e_{k-2}(u)+\dots +(-1)^kp_k(u)= 0, \quad k\ge 1.
\label{inv-Newt}
\ee
The highest order elementary symmetric polynomial (called {\em determinant})  commutes with the generating matrix $L(u)$
$$
e_m(u)L(v) = L(v)e_m(u),\quad \forall\,u,v,
$$
therefore, the Laurent coefficients of $e_m(u)$ are central elements of $\Y(R)$.
\end{proposition}
\begin{remark}\rm
If the involutive symmetry $R$ coincides with the flip $P$, then the braided Yangian $\Y(R)$ is nothing but the Drinfeld Yangian $\Y(gl(N))$  and
Proposition \ref{prop:inv-ident} reproduces the well-known results for it (see \cite{Mo}).
\end{remark}

\section{Bosonization of braided Yangians}

In this section we realize a bosonization of the braided Yangians, i.e. represent these Yangians in a braided analog of the Fock space. This
representation is also based on the evaluation map.

First, introduce the Fock space and  the corresponding bosonic algebra, associated with a skew-invertible Hecke symmetry $R$. The role of the Fock space ${\cal F}(R)$ is played by
the $R$-sym\-metric algebra ${\rm Sym}_R(V)$ of the basic space $V$. The elements of  the basis $\{x_i\}$ considered in section 2, generate this algebra. Also, the unit of this algebra
plays the role of the vacuum vector.

In the space $V^*$ we use the both bases: the right dual basis $\{x^j\}$ and the left dual one $\{\tx^j\}$, i.e. such that $<x^j, x_i>=\de_i^j$. The symmetric algebra  of the space $V^*$ is
defined by
\be
{\rm Sym}_R(V^*)=T(V^*)/\langle q x^i x^j-R^{ji}_{lk} x^k x^l \rangle.
\label{symdual}
\ee

Emphasize that in the left dual basis the algebra ${\rm Sym}_R(V^*)$ is defined by the same formula (up to replacing $x^i$ by $\tx^i$). To prove this claim we use the following relation
between two dual bases (see \cite{GPS})
\be
x^i=\tx^k B_k^i \label{bases}
\ee
and use the fact that  the product $B_1 B_2$  commutes with the braiding $R$.

Besides, we introduce the {\em permutation relations} between the generators of the algebras ${\rm Sym}_R(V^*)$ and ${\rm Sym}_R(V)$ by putting
\be
x^i x_j=B^i_j+\qq x_k x^l \Psi_{lj}^{ki},
\label{bos}
\ee
where $\Psi$ is defined in (\ref{def:Psi}).

\begin{proposition} The ideals coming in the definitions of the algebras ${\rm Sym}_R(V)$ and  ${\rm Sym}_R(V^*)$ are preserved by the permutation relations
{\rm (\ref{bos})}. Consequently, the map\rm
\be
{\rm Sym}_R(V^*)\ot {\rm Sym}_R(V)\to {\rm Sym}_R(V)\ot {\rm Sym}_R(V^*), \label{mappp}
\ee
\it
arising from successive permutations of all factors of a monomial $p\in  {\rm Sym}_R(V^*)$ and those of a monomial $q\in {\rm Sym}_R(V)$, is well-defined.
\end{proposition}

\noindent
{\bf Proof.}  Here we use  the left dual basis. By means of the map (\ref{bases}), we can write the relation (\ref{bos}) in terms of the left dual basis. Namely, we have
 \be
\tx^i x_j=\de_j^i+\qq x_k \tx^l (R^{-1})_{jl}^{ik}.
\label{bos1}
\ee

The fact that the systems (\ref{bos}) and (\ref{bos1}) are equivalent can be easily deduced from the relation (see \cite{O})
$$
B_1\,\Psi_{12}=R^{-1}_{21}\, B_2 \quad \Leftrightarrow \quad B_i^j \, \Psi_{jk}^{lm}=(R^{-1})_{ki}^{pl}\, B_p^m.
$$
With the use of (\ref{bos1}) and the Yang-Baxter equation for the Hecke symmetry $R$ we get:
$$
\tx^k(q\,x_i\,x_j-R_{ij}^{\,ab}x_a\,x_b) =q^{-2}(R^{-1})_{js}^{\,rb}(R^{-1})_{ir}^{\,ka}(q\,x_a\,x_b - R_{ab}^{\,cd}x_c\,x_d)\tx^s,
$$
This means that the permutation relations preserve the ideal generated by the relations
$$
q\,x_i\,x_j-R_{ij}^{\,ab}x_a\,x_b,
$$
therefore, the permutation relations are compatible with the algebraic structure of ${\rm Sym}_R(V)$. The same is true for
the ideal in the quotient (\ref{symdual}). Details are left to the reader. \hfill \rule{6.5pt}{6.5pt}

\medskip

Now, convert the generators of the algebra ${\rm Sym}_R(V^*)$ into operators acting on the space ${\cal F}(R)$.
Thus, in order to apply $\tx^i$  to an element $f\in {\cal F}(R)$ we send the product $\tx^i \ot f$ to the space
${\rm Sym}_R(V)\ot {\rm Sym}_R(V^*)$ by the map (\ref{mappp}). Then we apply the counit $\varepsilon$ to the factors
from ${\rm Sym}_R(V^*)$. This counit kills all elements of the algebra  ${\rm Sym}_R(V^*)$ except for those of
 the ground field $\K$ and acts as the identity operator on $\K$. Thus, we get an element  $g\in {\rm Sym}_R(V)\ot \K\cong {\rm Sym}_R(V)$.  Finally, we define
 the action as $\tx^i(f)=g$.   In particular, we have $\tx^j(1)=0$, $\tx^j(x_i)=\de_i^j$.

 Now,  define the creation and annihilation operators $ a_i^+,a^i \in{\rm End}({\cal F}(R))$ acting on the space ${\cal F}(R)$ as follows
\be
a_i^+ (f)= x_i\, f,\quad a^i(f)=\tx^i(f),\quad \forall\, f\in {\cal F}(R).
\label{anni}
\ee

\begin{proposition} The creation and annihilation operators meet the following permutation relations\rm
\be
q a_i^+ a_j^+=R_{ij}^{kl}a_k^+ a_l^+,\quad q a^i a^j=R_{lk}^{ji}a^k a^l,\quad a^j\,a_i^+-\qq\,  (R^{-1})_{ik}^{jl}  \, a_l^+\, a^k=\de_i^j.
\label{anni1}
\ee
\end{proposition}

\noindent
{\bf Proof.} The first relation is clear. The second relation follows from the fact that the permutation of  elements of the algebra ${\rm Sym}_R(V^*)$ with a fixed element
$f\in {\rm Sym}_R(V)$ is a morphism of the algebra ${\rm Sym}_R(V)$ onto itself. Also, the counit $\varepsilon$ is a morphism. So, the map $\tx^i \mapsto a^i$ is a representation
of the algebra ${\rm Sym}_R(V^*)$. The third relation is valid in virtue of the definition of the annihilation operators.  \hfill\rule{6.5pt}{6.5pt}

\medskip

The algebra generated by the creation and annihilation operators is called the {\em $R$-bosonic algebra} and is denoted ${\cal B}(R)$.

Now, we pass to the bosonization of the modified RE algebra $\lro$. 

\begin{proposition}
The elements $l_i^j=a_i^+ a^k B^j_k$ meet the defining relations of the algebra $\lro$. Consequently, the map
$$
\pi: l_i^j \mapsto a_i^+ a^k B^j_k,
$$
which is called a bosonization of the algebra $\lro$, defines a representation of this algebra in the Fock space ${\cal F}(R)$.
\end{proposition}

\noindent
{\bf Proof.} Here, we use the notations $x_i$ and $x^j$ instead of $a_i^+$ and $a^p B^j_p$ respectively.
Let us present the relations (\ref{bos}) in an equivalent form
$$
x^a R_{ai}^{bj}x_b=\de_i^j+\qq x_i x^j.
$$

Now, transform the matrix $\pi(R_1L_1R_1L_1)$ as follows (underlined are the terms affected by transformations on the next step)
\begin{eqnarray*}
\pi(R_{ij}^{ab}l_a^cR_{cb}^{dl}l_d^{k}) &=&R_{ij}^{ab}x_a\underline{x^c R_{cb}^{dl}x_d}\,x^{k} =
R_{ij}^{al}x_ax^{k}+q^{-1}\underline{R_{ij}^{ab}x_a\,x_b}\,x^{l}x^{k} =\\
 &=& \pi(R_{ij}^{al}l_a^{k}) + x_{i}x_{j}x^{l}x^{k}.
\end{eqnarray*}
Analogously, we transform the entries of the matrix $L_1R_1L_1R_1$:
\begin{eqnarray*}
\pi(l_{i}^aR_{aj}^{bc}l_b^dR_{dc}^{k l})&=&x_{i}\underline{x^aR_{aj}^{bc}x_b}x^dR_{dc}^{k l}=
x_{i}x^dR_{dj}^{k l}+ q^{-1}x_{i}x_{j}\underline{x^cx^dR_{dc}^{k l}}=\\
 &=& \pi(l_{i}^dR_{dj}^{k l}) + x_{i}x_{j}x^{l}x^{k}.
\end{eqnarray*}
Subtracting the latter expression from the former one we come to the result. \hfill \rule{6.5pt}{6.5pt}

\medskip

Now, we are able to realize a bosonization of the braided Yangians based on the evaluation morphism. Let us represent the generating matrix $L(u)$ of a given braided
Yangian $\YR$ via the map $L(u)\mapsto I+\frac{M}{u}$, where $M$ is either the generating matrix of the RE algebra $\lr$, if $R$ is a Hecke symmetry, or that of $\lro$ if
$R$ is an involutive symmetry.

First, consider the latter case. Using the bosonization of the algebra $\lro$ described above, we immediately get a bosonization of the braided Yangian $\YR$. More precisely, we have the
following proposition.

\begin{proposition}
The map\rm
\be
l_i^j(u)\mapsto \de_i^j+ \frac{1}{u}\,a_i^+a^k B_k^j,
\label{map2}
\ee
\it
where $a_i^+$ and $a^j$ are respectively the creation and annihilation operators acting on the Fock space ${\cal F}(R)$
and subject to {\rm (\ref{anni1})} with $q=1$, is a representation of the braided Yangian $\Y(R)$, associated with an involutive symmetry $R$.
\end{proposition}

In terms of the Laurent coefficients  the map (\ref{map2}) can be presented as follows
$$
l_i^j[1]\mapsto a_i^+ a^k B_k^j,\quad  L[k]\mapsto 0,\,\, k\geq 2.
$$

Note that  the denominator in (\ref{map2}) can be replaced by $u-u_0$, $u_0\in \K$. Such representations
with different $u_0$ can be multiplied by means of the coproduct, defined in the braided Yangians. Note that the vacuum vector of the product of the corresponding Fock
spaces is the product of the vacuum vectors of the factors.

If $R$ is a Hecke symmetry, then by using the map inverse to (\ref{ma}) with $h=1$, we arrive to the following bosonization of the braided Yangian $\Y(R)$
(also, presented by means  of the Laurent coefficients)
\be
l_i^j[1]\mapsto a_i^+a^k B_k^j+\frac{\de_i^j}{q-\qq},\quad  L[k]\mapsto 0,\,\, k\geq 2.
\label{mapss}
\ee

\begin{proposition}
The map {\rm (\ref{mapss})} is a representation of the braided Yangian $\YR$ associated with a Hecke symmetry $R$.
\end{proposition}

This map does not admit a limit $q\to 1$ since the evaluation maps for the braided Yangians corresponding
to involutive symmetries and Hecke ones essentially differ. However, it is possible by renormalizing  the map (\ref{mapss}) to get the following representation
$$ 
l_i^j[1]\mapsto (q-\qq) a_i^+a^k B_k^j+\de_i^j,\quad L[k]\mapsto 0,\,\, k\geq 2,
$$
which has a limit. However, this limit is out of interest.

In the classical case, i.e. while $R(u,v)$ is the Yang $R$-matrix, the role of the bosonic algebra is often attributed to the  Zamolodchikov-Faddeev (ZF) algebra
(see \cite{ZZ, F}). It appeared   in the framework of  the second quantization of some dynamical models. In particular, it comes to the Nonlinear Schrodinger hierarchy (see \cite{MRSZ}).
Nevertheless, a consequent definition of a bosonic algebra in the spirit of the ZF one meets some difficulties, which we comment below. Also, we suggest a modified version of the
bosonic algebra valid for all current $R$-matrices in question.

As was pointed out in \cite{LM}, a consistent Fock space can be only associated with an involutive current $R$-matrix, i.e. that subject to the condition: $R(u,v)\, R(v,u)=I$.
The authors of \cite{LM} found such a current $R$-matrix in the form $g(u,v)I$, where $g(u,v)$ is an appropriate function. However, this braiding is out of our interest.

So, assuming $R(u,v)$ to be one of the  current $R$-matrices constructed in section 4, we first normalize it  as suggested there, and get an involutive braiding
$\RRR(u,v)$, given by one of formulae (\ref{demo}). Then we define the Fock space  ${\cal F}(\RRR)$ as an $\RRR$-symmetric algebra of a space $V(w)$.
The space ${\cal F}(\RRR)$ consists of finite linear combinations of products  $x_{i_1}(u_1)\ot x_{i_2}(u_2)\ot\dots \ot  x_{i_k}(u_k)$, modulo the ideal generated
by quadratic elements
\be
x_{i}(u)x_{j}(v)=\RRR_{ij}^{kl}(u,v) \,x_{k}(v) x_{l}(u).
\label{ZF1}
\ee

Rewrite these relations in a more detailed form
\be
q\,x_{i}(u)x_{j}(v)-R_{ij}^{kl}\, x_{k}(v) x_{l}(u)=g(u,v)(x_{i}(u) x_{j}(v)- x_{i}(v) x_{j}(u))
\label{ZF0}
\ee
where $g(u,v)$ is one of the  functions
$$
g(u,v)=\frac{1}{u-v}\qquad{\rm or}\qquad g(u,v)=\frac{u(q-\qq)}{u-v}.
$$
In the former case we also set $q=1$ in (\ref{ZF0}).

In a similar manner, we define the $\RRR$-symmetric algebra of the dual space $V^*(w)$ by imposing the following system of relations:
\be
\tx^l(u)\tx^k(v)=\RRR_{ij}^{kl}(u,v)\, \tx^j(v)\tx^i(u).
\label{ZF2}
\ee

Now, according to the above pattern, define the corresponding permutation relations by the  formula similar to (\ref{bos}):
\be
\tx^j(u) x_i(v)=\RRR_{ik}^{jl}(v,u)\,x_l(v) \tx^k(u)+\de_i^j \de(u-v).
\label{ZF3}
\ee

We leave to the reader checking that these permutation relations preserve the ideals coming in the definitions of the above symmetric algebras.
Finally, we define the corresponding bosonic algebra by the relations (\ref{ZF1}), (\ref{ZF2}), and (\ref{ZF3}).

Observe that
even in the case when $R(u,v)$ is the Yang $R$-matrix, our definition of the bosonic algebra differs  from that of the ZF algebra.

\bigskip
\noindent
{\bf Acknowledgement.} D.G. is grateful to
Indiana University – Purdue University Indianapolis for
stimulating atmosphere during his one month scientific stay at Department of Mathematical Sciences in 2015 year, this paper has been conceived then. He is also
thankful to Professors Vitaly Tarasov and Alexander Molev for elucidating discussions.

The work of P.S. has been funded by the Russian Academic Excellence Project '5-100' and by the RFBR grant no. 16-01-00562-A.

\end{document}